\documentclass{article}
\usepackage{graphicx}
\usepackage{amscd}
\usepackage[matrix,arrow,curve]{xy}
\usepackage{amssymb,amsmath}
\usepackage{amsmath,amsfonts,amssymb,latexsym}
\newtheorem{theo}{Theorem}[section]
\newtheorem{mydf}{Definition}[section]
\newtheorem{prop}{Proposition}[section]
\newtheorem{cor}{Corollary}[section]
\newtheorem{lem}{Lemma}[section]

\newtheorem{exam}{Example}[section]

\def\P{\bf Proof. \rm}

\newcommand{\proofterminator}{\hfill\ensuremath{\square}}
\begin{document}
\title{\Large{Fundamentals for Symplectic $\mathcal{A}$-modules}}
\author{{\small \underline{Anastasios Mallios}} \\ {\small  Department of Mathematics} \\
{\small University of Athens}\\ {\small Athens, Greece}\\
{\small Email: amallios@math.uoa.gr}\\ {\small \underline{Patrice P. Ntumba}}\\
{\small Department of Mathematics and Applied
Mathematics}\\{\small University of Pretoria}\\ {\small Hatfield
0002, Republic of South Africa}\\{\small Email:
patrice.ntumba@up.ac.za}\\}

\date{}
\maketitle

\begin{abstract} In his \cite{mallios}, \cite{mallios1}, the first
author shows that the \textit{sheaf theoretically} based
\textit{Abstract Differential Geometry incorporates and generalizes
all the classical differential geometry}. Here, we undertake to
partially explore the implications of Abstract Differential Geometry
to classical \textit{symplectic geometry}. The full investigation
will be presented elsewhere.
\end{abstract}

\maketitle

\section*{Introduction} Here is an attempt of taking the theory of
\textit{Abstract Differential Geometry} ($\grave{a}$ la Mallios),
\cite{mallios}, \cite{mallios1}, to new horizons, such as those
related to the classical \textit{symplectic geometry}. In this
endeavor, we show the \textit{extent} to which tools provided by
Abstract Differential Geometry help (re)capture in a sheaf-theoretic
manner fundamental notions and results which characterize the
standard \textit{symplectic algebra}. This endeavor will pave the
way to rewrite and/or recover a great deal of classical symplectic
geometry, with no use at all of any notion of ``differentiability"
(differentiability is here understood in the sense of the standard
\textit{differential geometry} of $C^\infty$-manifolds). As a result
of all \textit{this}, we show that \textit{sheaf-theory methods}
turn out to be an appropriate way for the \textit{algebraization} of
classical symplectic geometry. As pointed out by the first author in
\cite{mallios1}, \cite{mallios2}, this algebraic approach to
differential (and symplectic) geometry is of a particular interest
to theoretical physicists, for it has been for long the demand
and/or wish of many to ``\textit{find a purely algebraic theory for
the description of reality}" (i.e. physics according to our
understanding) (A. Einstein, [\cite{einstein}, p.166]).

Our main reference, throughout the present account, is the first
author's book \cite{mallios}, for which the reader is requested to
have handy for we have skipped some necessary \textit{basics} of
Abstract Differential Geometry.

The paper is divided into four sections. \S 1 concerns with
$\mathcal{A}$-\textit{tensors}: these are the counterparts, in this
framework, of \textit{classical} tensors. Some results, pertaining
to the standard \textit{multilinear algebra machinery}, are hereby
provided. $\mathcal{A}$-tensors constitute a precursor to the
fundamental theory of \textit{exterior}
$\mathcal{A}$-$k$-\textit{forms}, which are developed in \S 2. The
\textit{exterior algebra sheaf} is defined as the \textit{direct sum
of sheaves of germs of exterior $\mathcal{A}$-$k$-forms on an
$\mathcal{A}$-module $\mathcal{E}$}; this sum is endowed with the
exterior product $\wedge$. In \S 3, we show that given a
\textit{non-zero skew-symmetric non-degenerate
$\mathcal{A}$-morphism} $\omega$ on the \textit{standard free
$\mathcal{A}$-module of rank} $n$, defined on a topological space
$X$, there is a \textit{basis} $\mathcal{B}$ of $\mathcal{A}^n(X)$,
relative to which the matrix of $\omega$ is
\[\left[\begin{array}{cc} 0 & \mbox{I}_n \\ -\mbox{I}_n &
0\end{array}\right],\]where $\mbox{I}_n$ is the $n\times n$ identity
matrix (the diagonal entry $1$ in $I$ is the \textit{global identity
section}). It also turns out that this result holds for
\textit{vector sheaves} as well. We further introduce in this same
section the notion of \textit{symplectic group sheaf} of an
arbitrary $\mathcal{A}$-module. In \S 4, we deal with
\textit{characteristic polynomial, eigenvector and eigenvalue
sections} of an $\mathcal{A}$-module, and prove the corresponding
version of the \textit{Cayley-Hamilton} theorem.

\section{$\mathcal{A}$-Tensors} Throughout this paper, the pair $(X,
\mathcal{A})$, or just $\mathcal{A}$ will denote a fixed
$\mathbb{C}$-\textit{algebraized space}, where $X$ is a
topological space and $\mathcal{A}$ a \textit{sheaf} (over $X$) of
\textit{unital, commutative algebras}. For more details about
algebraized spaces, see \cite{mallios}.

Let \[\mathcal{E}\equiv (\mathcal{E}, \pi, X)\]be an
$\mathcal{A}$-module on $X$. The (complete) \textit{presheaves of
sections} of sheaves $\mathcal{A}$ and $\mathcal{E}$ are denoted
by
\[\begin{array}{lll} \Gamma(\mathcal{A})\equiv (\Gamma(U,
\mathcal{A}), \tau^U_V) & \mbox{and} & \Gamma(\mathcal{E})\equiv
(\Gamma(U, \mathcal{E}), \pi^U_V) ,\end{array}\]respectively. Let
$\mathcal{E}^\ast$ be the dual $\mathcal{A}$-module of
$\mathcal{E}$; so \[\mathcal{E}^\ast=
\mathcal{H}om_\mathcal{A}(\mathcal{E}, \mathcal{A})\](see
\cite{mallios}, p. 298). It is easy (cf. \cite{mallios}, p. 129) to
see that the correspondence that associates with every open subset
$U$ of $X$ the $\Gamma(U, \mathcal{A})$-module
\[\otimes^p\Gamma(U, \mathcal{E})\otimes_{\Gamma(U,
\mathcal{A})}\otimes^q\Gamma(U, \mathcal{E}^\ast),\]where
\[\otimes^p\Gamma(U, \mathcal{E}):= \underbrace{\Gamma(U, \mathcal{E})\otimes_{\Gamma(U,
\mathcal{A})}\ldots \otimes_{\Gamma(U, \mathcal{A})}\Gamma(U,
\mathcal{E})}_{\mbox{$p-$times}},\] and
\[\otimes^q\Gamma(U, \mathcal{E}^\ast):= \underbrace{\Gamma(U,
\mathcal{E}^\ast)\otimes_{\Gamma(U, \mathcal{A})}\ldots
\otimes_{\Gamma(U, \mathcal{A})}\Gamma(U,
\mathcal{E}^\ast)}_{\mbox{$q-$times}},\]along with the obvious
restriction morphisms, provides a $\textit{presheaf of
$\Gamma(\mathcal{A})$-modules}$ on $X$. We denote this presheaf by
\begin{equation}\label{tensor}T^p_q\Gamma(\mathcal{E})\equiv
\otimes^p\Gamma(\mathcal{E})\otimes_{\Gamma(\mathcal{A})}\otimes^q\Gamma(\mathcal{E}^\ast).\end{equation}

\begin{mydf}\label{shtensor} \emph{Given an $\mathcal{A}$-\textit{module} $\mathcal{E}$ on a
topological space $X$, we denote by
\[\otimes^p\mathcal{E}\otimes_\mathcal{A}\otimes^q\mathcal{E}^\ast\]
the \textit{sheaf} generated by the presheaf
$T^p_q\Gamma(\mathcal{E})$, given in $(\ref{tensor})$, i.e. one
has
\[\otimes^p\mathcal{E}\otimes_\mathcal{A}\otimes^q\mathcal{E}^\ast:=
\mathbf{S}(T^p_q\Gamma(\mathcal{E}))=
\mathbf{S}(\otimes^p\Gamma(\mathcal{E})\otimes_{\Gamma(\mathcal{A})}\otimes^q\Gamma(\mathcal{E}^\ast)).\]}\end{mydf}

We extend the definition of
$\otimes^p\mathcal{E}\otimes_\mathcal{A}\otimes^q\mathcal{E}^\ast$
to the cases $p=0$ and $q=0$ by setting
\[\otimes^0\mathcal{E}\otimes_\mathcal{A}\otimes^q\mathcal{E}^\ast=
\otimes^q\mathcal{E}^\ast,\]and
\[\otimes^p\mathcal{E}\otimes_\mathcal{A}\otimes^0\mathcal{E}^\ast=
\otimes^p\mathcal{E}.\]The elements of
$\otimes^p\mathcal{E}\otimes_\mathcal{A}\otimes^q\mathcal{E}^\ast$
are called $\textbf{$\mathcal{A}$-tensors}$ over $\mathcal{E}$,
and are said to be $\textit{contravariant of order $p$}$ and
$\textit{covariant of order $q$}$; or simply, of $\textit{type
$(p, q)$}$.

The next lemma shows the analog of a classical result of (ordinary)
tensors of type $(p, q)$. Before we examine the result, let us
recall that given an algebra sheaf $\mathcal{A}$ on a topological
$X$, by a \textit{vector sheaf $\mathcal{E}$, of a rank $n$}, on
$X$, we mean a \textit{locally free $\mathcal{A}$-module of rank
$n$} on $X$; that is for every $x\in X$, there exists an open
neighborhood $U$ of $x\in X$ such that one has
\begin{equation}\label{gauge}\mathcal{E}|_U= \mathcal{A}^n|_U,\end{equation} with the equality
sign being actually an $\mathcal{A}|_U$-isomorphism of the
$\mathcal{A}|_U$-modules $\mathcal{E}|_U$ and $\mathcal{A}^n|_U.$
Any open set $U$ in $X$ for which (\ref{gauge}) holds is called a
local gauge of $\mathcal{E}$.

\begin{lem}\label{mult} \emph{If $\mathcal{E}$ is a \textit{vector sheaf} on a topological
space $X$, then
\[\otimes^p\mathcal{E}\otimes_\mathcal{A}\otimes^q\mathcal{E}^\ast=
\mathcal{H}om_\mathcal{A}(\otimes^p\mathcal{E}^\ast\otimes_\mathcal{A}\otimes^q\mathcal{E},
\mathcal{A}).\]}\end{lem} \P This is an easy verification. In
fact, based on Mallios[\cite{mallios}, Comment (5.27), p. 132,
Theorems 5.1, p. 299, 6.1, p. 302, 6.2, p. 304, and Corollary 6.2,
p. 305], one has
\begin{eqnarray*}\otimes^p\mathcal{E}\otimes_\mathcal{A}\otimes^q\mathcal{E}^\ast
& = &
(\otimes^p\mathcal{E})^{\ast\ast}\otimes_\mathcal{A}\otimes^q\mathcal{E}^\ast
\\ & = & \mathcal{H}om_\mathcal{A}((\otimes^p\mathcal{E})^\ast,
\otimes^q\mathcal{E}^\ast)\\ & = &
\mathcal{H}om_\mathcal{A}(\otimes^p\mathcal{E}^\ast,
(\otimes^q\mathcal{E})^\ast) \\ & = &
\mathcal{H}om_\mathcal{A}(\otimes^p\mathcal{E}^\ast,
\mathcal{H}om_\mathcal{A}(\otimes^q\mathcal{E}, \mathcal{A}))\\ &
= &
\mathcal{H}om_\mathcal{A}(\otimes^p\mathcal{E}^\ast\otimes_\mathcal{A}\otimes^q\mathcal{E},
\mathcal{A}).\end{eqnarray*}\proofterminator

A corollary that one can derive from Lemma \ref{mult} requires the
following definitions.

\begin{mydf} \emph{Let $\mathcal{E}_1,\ldots, \mathcal{E}_n$, $n\in
\mathbb{N}$, and $\mathcal{F}$ be $\mathcal{A}$-\textit{modules}
on the same topological space $X$. The $\mathcal{A}$-morphism
$\varphi: \mathcal{E}_1\times\ldots \times
\mathcal{E}_n\longrightarrow \mathcal{F}$ is called an
\textbf{$\mathcal{A}$-multilinear morphism} if, for all
\textit{open subset} $U\subseteq X$,
\[\varphi_U: \Gamma(U, \mathcal{E}_1)\times_{\Gamma(U, \mathcal{A})}\ldots \times_{\Gamma(U, \mathcal{A})}
\Gamma(U, \mathcal{E}_n)\longrightarrow \Gamma(U, \mathcal{F})\]is a
$\Gamma(U, \mathcal{A})$-\textit{multilinear morphism} for the
$\Gamma(U, \mathcal{A})$-\textit{modules concerned}.}\end{mydf} We
are now set to \textit{generalize the functor}
$\mathcal{H}om_\mathcal{A}$ ($\mathcal{H}om_\mathcal{A}$ is a
bifunctor $\mathcal{A}$-$\mathcal{M}od_X\longrightarrow
\mathcal{A}$-$\mathcal{M}od_X$, where
$\mathcal{A}$-$\mathcal{M}od_X$ is the category of
$\mathcal{A}$-modules on a topological space $X$, (see
\cite{mallios}, p. 133)), to \textit{the functor}
$\mathcal{L}_\mathcal{A}^n$, $n\in\mathbb{N}$, which we define
below.

In effect, suppose that we are given $\mathcal{A}$-modules
$\mathcal{E}_i$, $i=1, \ldots, n$, and $\mathcal{F}$ on the same
topological space $X$. For any open set $U$ in $X$, let
\begin{equation}\label{Hom}\mbox{Hom}_{\mathcal{A}|_U}({\mathcal{E}_1}|_U\times
\ldots\times{\mathcal{E}_n}|_U, \mathcal{F}|_U)\equiv
L^n_{\mathcal{A}|_U}({\mathcal{E}_1}|_U,
\ldots,{\mathcal{E}_n}|_U; \mathcal{F}|_U)\end{equation}be the set
of $\textit{$\mathcal{A}|_U$-$n$-linear morphisms}$ of the
$\mathcal{A}|_U$-module ${\mathcal{E}_1}\times
\ldots\times{\mathcal{E}_n}|_U$ into the $\mathcal{A}|_U$-module
$\mathcal{F}|_U$.

\begin{lem}\emph{The set, in $(\ref{Hom})$, is a \textit{module} over $\mathcal{A}(U)$; hence, in
particular, a $\mathbb{C}$-\textit{vector space}.} \end{lem} \P
The proof is similar to the proof of Statement 6.1, p. 133,
\cite{mallios}. \proofterminator

On the other hand, it is readily verified that, given
$\mathcal{A}$-modules $\mathcal{E}_1,\ldots, \mathcal{E}_n$, and
$\mathcal{F}$ on $X$ as above, the correspondence
\begin{equation}\label{Hom1}U\longmapsto L^n_{\mathcal{A}|_U}({\mathcal{E}_1}|_U,\ldots,
{\mathcal{E}_n}|_U; \mathcal{F}|_U),\end{equation}where $U$ runs
over the open subsets of $X$, along with the obvious
\textit{restriction maps} yields a \textit{complete presheaf of
$\mathcal{A}$-modules} on $X$.

Thus, we have

\begin{mydf} \emph{Let $\mathcal{E}_1, \ldots, \mathcal{E}_n$ and
$\mathcal{F}$ be $\mathcal{A}$-\textit{modules} on a topological
space $X$. By the $\textbf{sheaf of germs of
$\mathcal{A}$-$n$-linear morphisms}$ of $\mathcal{E}_1\times
\ldots\times \mathcal{E}_n$ in $\mathcal{F}$, we mean the
\textit{sheaf}, on $X$, generated by the (complete)
\textit{presheaf}, defined by $(\ref{Hom1})$. We denote the
induced sheaf by
\[\mathcal{L}_\mathcal{A}^n(\mathcal{E}_1,\ldots, \mathcal{E}_n;
\mathcal{F}).\]}\end{mydf}

We may now state

\begin{cor}\emph{Let $\mathcal{E}$ be an $\mathcal{A}$-\textit{module }on $X$.
Then,
\[\otimes^p\mathcal{E}\otimes_\mathcal{A}\otimes^q\mathcal{E}^\ast =
\mathcal{L}_\mathcal{A}^{p+q}(\underbrace{\mathcal{E}^\ast,\ldots,
\mathcal{E}^\ast}_{\mbox{$p$-times}},
\underbrace{\mathcal{E},\ldots, \mathcal{E}}_{\mbox{$q$-times}};
\mathcal{A})\equiv \mathfrak{T}^p_q(\mathcal{E}),\]where
$\mathcal{L}_\mathcal{A}^{p+q}(\mathcal{E}^\ast,\ldots,
\mathcal{E}^\ast, \mathcal{E},\ldots, \mathcal{E}; \mathcal{A})$
is the $\mathcal{A}$-\textit{module of
$\mathcal{A}$-$(p+q)$-linear morphisms}.}\end{cor} \P Using Lemma
\ref{mult} and Mallios[\cite{mallios}, Lemma 5.1, p. 132 and
Definition 6.1, p. 134], one has, for every open set $U$ in $X$,
\[\begin{array}{ll}
\mathcal{H}om_\mathcal{A}(\otimes^p\mathcal{E}^\ast\otimes_\mathcal{A}\otimes^q\mathcal{E},
\mathcal{A})(U) & =
\mbox{Hom}_{\mathcal{A}|_U}((\otimes^p\mathcal{E}^\ast\otimes_\mathcal{A}\otimes^q\mathcal{E})|_U,
\mathcal{A}|_U) \\ &
=\mbox{Hom}_{\mathcal{A}|_U}(\otimes^p(\mathcal{E}^\ast|_U)\otimes_{\mathcal{A}|_U}\otimes^q(\mathcal{E}|_U),
\mathcal{A}|_U)\\ & =L^{p+q}_{\mathcal{A}|_U}(\mathcal{E}^\ast|_U,
\ldots, \mathcal{E}^\ast|_U, \mathcal{E}|_U, \ldots, \mathcal{E}|_U;
\mathcal{A}|_U).\end{array}\]Therefore,\[\mathcal{H}om_\mathcal{A}(\otimes^p\mathcal{E}^\ast
\otimes_\mathcal{A}\otimes^q\mathcal{E}, \mathcal{A}) =
\mathcal{L}_\mathcal{A}^{p+q}(\mathcal{E}^\ast,\ldots,
\mathcal{E}^\ast, \mathcal{E},\ldots, \mathcal{E};
\mathcal{A}).\]\proofterminator

So, we come now to the following definition.

\begin{mydf} \emph{Let $\mathcal{E}$ be an $\mathcal{A}$-\textit{module} on a
topological space $X$, and let $\mathbf{t}_1\in
\mathfrak{T}^{i_1}_{j_1}(\mathcal{E})$ and $\mathbf{t}_2\in
\mathfrak{T}^{i_2}_{j_2}(\mathcal{E})$. The
$\mathcal{A}$-$\textbf{tensor product}$ of $\mathbf{t}_1$ and
$\mathbf{t}_2$ is the $\mathcal{A}$-\textit{tensor}
$\mathbf{t}_1\otimes \mathbf{t}_2\in
\mathfrak{T}^{i_1+i_2}_{j_1+j_2}(\mathcal{E})$, defined by
\begin{eqnarray*}\lefteqn{\mathbf{t}_1\otimes \mathbf{t}_2(s_1, \ldots, s_{i_1}, t_1,
\ldots, t_{i_2}, u_1, \ldots, u_{j_1}, v_1, \ldots, v_{j_2})} \nonumber\\
& & = \mathbf{t}_1(s_1, \ldots, s_{i_1}, u_1, \ldots,
u_{j_1})\mathbf{t}_2(t_1, \ldots, t_{i_2}, v_1, \ldots,
v_{j_2})\end{eqnarray*}where $s_\alpha, t_\alpha\in
\mathcal{E}^\ast(U)$ and $u_\beta, v_\beta\in \mathcal{E}(U)$, and
where, for all $k=1,2$, the $\mathcal{A}-$tensor $\mathbf{t}_k$,
viewed as a map on \textit{sections} of the $\mathcal{A}$-modules of
\[\begin{array}{ll}
\underbrace{\mathcal{E}^\ast\times_\mathcal{A}\ldots\times_\mathcal{A}
\mathcal{E}^\ast}_{i_k}\times_\mathcal{A}
\underbrace{\mathcal{E}\times_\mathcal{A}\ldots\times_\mathcal{A}\mathcal{E}}_{j_k}
& {\mbox{and $\mathcal{A},$}}\end{array}\]is the
$\mathcal{A}(U)$-$(i_k+j_k)$-linear morphism
\[
\mathcal{E}^\ast(U)\times_{\mathcal{A}(U)}\ldots\times_{\mathcal{A}(U)}\mathcal{E}^\ast(U)
\times_{\mathcal{A}(U)}\mathcal{E}(U)
\times_{\mathcal{A}(U)}\ldots
\times_{\mathcal{A}(U)}\mathcal{E}(U)\longrightarrow
\mathcal{A}(U),\]for all \textit{open} subset $U\subseteq X$.}
\end{mydf}

One can be assured that the standard $\textit{multilinear algebra
machinery}$ can be appropriately reformulated within the present
setting. For example, let us look at

\begin{prop}\emph{Let $\mathcal{E}$ be a \textit{vector sheaf}
of rank $n$ on a topological space $X$. Then, for all $k, l\in
\mathbb{N}$, the $\mathcal{A}$-module
$\mathfrak{T}^k_l(\mathcal{E})$ is a vector sheaf of rank
$n^{k+l}$.} \end{prop} \P The proof is based on relations (5.25), p.
132, (6.23), p. 137, and Statement (5.19), p. 301, all found in
\cite{mallios}.

Let $x\in X$ and $U$ an \textit{open} neighborhood of $x$ such
that $\mathcal{E}|_U= \mathcal{A}^n|_U= \mathcal{E}^\ast|_U$.
Then, for all open subset $V\subseteq U$, one has
\begin{eqnarray*}(\mathfrak{T}^k_l(\mathcal{E})|_U)(V) & = &
(\mathcal{L}_\mathcal{A}^{k+l}(\underbrace{\mathcal{E}^\ast, \ldots,
\mathcal{E}^\ast}_k, \underbrace{\mathcal{E},\ldots, \mathcal{E}}_l;
\mathcal{A})|_U)(V)\\ & = &
\mathcal{L}_\mathcal{A}^{k+l}(\mathcal{E}^\ast, \ldots,
\mathcal{E}^\ast, \mathcal{E},\ldots, \mathcal{E}; \mathcal{A})(V)
\\ & = & L^{k+l}_{\mathcal{A}|V}(\mathcal{E}^\ast|_V, \ldots,
\mathcal{E}^\ast|_V, \mathcal{E}|_V, \ldots, \mathcal{E}|_V;
\mathcal{A}|_V) \\ & = &
\mbox{Hom}_{\mathcal{A}|_V}(\mathcal{E}^\ast\times\ldots\times\mathcal{E}^\ast|_V\times
\mathcal{E}|_V\times\ldots\times \mathcal{E}|_V, \mathcal{A}|_V)\\
& = &
\mbox{Hom}_{\mathcal{A}|_V}(\mathcal{A}^n|_V\times\ldots\times
\mathcal{A}^n|_V\times \mathcal{A}^n|_V\times \ldots\times
\mathcal{A}^n|_V, \mathcal{A}|_V) \\ & = &
\mbox{Hom}_{\mathcal{A}|_V}((\mathcal{A}^n\otimes\ldots\otimes
\mathcal{A}^n)|_V, \mathcal{A}|_V) \\ & = &
\mathcal{H}om_\mathcal{A}(\mathcal{A}^n\otimes
\ldots\otimes\mathcal{A}^n, \mathcal{A})(V)\\ & = &
\mathcal{H}om_\mathcal{A}({\mathcal{A}}^{n^{k+l}},
\mathcal{A})(V)\\ & = &
(\mathcal{A}^{n^{k+l}}|_U)(V),\end{eqnarray*}which shows that
\[\mathfrak{T}^k_l(\mathcal{E})|_U=
\mathcal{A}^{n^{k+l}}|_U,\]that is the $\mathcal{A}$-module
$\mathfrak{T}^k_l(\mathcal{E})$ is a vector sheaf of rank
$n^{k+l}$, as desired.\proofterminator

\begin{cor} \emph{Let $\mathcal{E}$ be a \textit{vector sheaf
} of rank $n$ on a topological space $X$, and $\{s_i\}_{1\leq i\leq
n}$ a \textit{basis} of the $\mathcal{A}(U)$-module $\Gamma(U,
\mathcal{E})$, with $U$ an open subset of $X$ such that
$\mathcal{E}|_U= \mathcal{A}^n|_U$. Then, for all $k, l\in
\mathbb{N}$, a basis of the $\mathcal{A}(U)$-\textit{module}
$\Gamma(U, \mathfrak{T}^k_l(\mathcal{E}))$ is given by
\[\{s_{i_1}\otimes \ldots\otimes s_{i_k}\otimes
{s^\ast}^{j_1}\otimes\ldots\otimes {s^\ast}^{j_l}|\ i_p,
j_p=1,\ldots, n\},\]where $\{{s^\ast}^j\}_{1\leq j\leq n}$ is the
dual basis of $\{s_i\}_{1\leq i\leq n}$.}\end{cor}\P The proof is
similar to the proof of Proposition 1.7.2, p. 53,
\cite{ralph}.\proofterminator

We close this section with the following important definition,
which will be of use in the sequel. (See also [\cite{mallios}: p.
301, (5.22)- (5.24)].)

\begin{mydf} \emph{Let $\mathcal{E}$ and $\mathcal{F}$ be
$\mathcal{A}$-\textit{modules} on a topological space $X$. For
$\varphi\in \mathcal{H}om_\mathcal{A}(\mathcal{E}, \mathcal{F})$,
in keeping with the classical notation, see $\cite{bourbaki}$, p.
$234$, or $\cite{chambadal}$, p. $68$, we define the
$\mathbf{transpose\ of}$ $\varphi$ by
\[{}^t\varphi\in \mathcal{H}om_\mathcal{A}(\mathcal{F}^\ast,
\mathcal{E}^\ast)\]such that
\begin{equation}\begin{array}{ll}({}^t\varphi)(u):=
u\circ\varphi, & u\in
\mathcal{F}^\ast(U)\end{array}\label{transpose}\end{equation}with
$U$ \textit{open} in $X$, i.e., in other words,\[\begin{array}{ll}
{}^t\varphi(u)(v)= u(\varphi(v)), & v\in
\mathcal{E}(U).\end{array}\] $($ In $(\ref{transpose})$, we have
used the ${}^t\varphi$-\textit{corresponding map on sections} of
the $\mathcal{A}$-modules $\mathcal{E}^\ast$ and
$\mathcal{F}^\ast$. $)$}\end{mydf}

\section{Exterior $\mathcal{A}$-$k$-forms} As has been the case so
far, we assume in this section as well that the triple
$(\mathcal{A}, \tau, X)$ stands for the \textit{sheaf of
commutative $\mathbb{C}$-algebras with an identity element} on a
topological space $X$. Furthermore, we let
\[\Gamma(\mathcal{A})\equiv (\Gamma(U, \mathcal{A}), \tau^U_V)\]be
the corresponding (complete) \textit{presheaf of sections} of
$\mathcal{A}$.

Now, let $\mathcal{E}$ be an $\mathcal{A}$-module on $X$. For any
\textit{open} set $U\subseteq X$, let
\[\Omega_{\mathcal{A}|_U}^k(\underbrace{\mathcal{E}|_U\oplus\ldots
\oplus\mathcal{E}|_U}_{k-\mbox{times}}, \mathcal{A}|_U)\]be the
\textit{set of all skew-symmetric $\mathcal{A}|_U$-$k$-linear
morphisms} of the $\mathcal{A}|_U$-modules
$\mathcal{E}|_U\oplus\ldots\oplus \mathcal{E}|_U$ and
$\mathcal{A}|_U$. It is obvious that we have
\[\Omega^k_{\mathcal{A}|_U}(\mathcal{E}|_U\oplus\ldots\oplus\mathcal{E}|_U,
\mathcal{A}|_U)\subseteq
\mbox{Hom}_{\mathcal{A}|_U}(\mathcal{E}|_U\oplus\ldots\oplus\mathcal{E}|_U,
\mathcal{A}|_U),\]where for every open set $U\subseteq X$,
\[\mbox{Hom}_{\mathcal{A}|_U}(\mathcal{E}|_U\oplus\ldots\oplus\mathcal{E}|_U,
\mathcal{A}|_U)\]is the set of all
$\mathcal{A}|_U$-$k$-\textit{linear morphisms} of
$\mathcal{E}|_U\oplus\ldots\oplus\mathcal{E}|_U$ into
$\mathcal{A}|_U$.

As is naturally expected, the correspondence
\begin{equation}\label{ext}U\longmapsto
\Omega^k_{\mathcal{A}|_U}(\mathcal{E}|_U\oplus\ldots\oplus\mathcal{E}|_U,
\mathcal{A}|_U),\end{equation}where $U$ is open in $X$, along with
the obvious \textit{restriction maps}, yields a \textit{complete
presheaf of $\mathcal{A}$-modules} on $X$. The sheaf, on $X$,
generated by the presheaf defined by $(\ref{ext})$ is called the
\textbf{sheaf of germs of exterior $\mathcal{A}$-$k$-forms on
$\mathcal{E}$}, and is denoted \[\Omega^k(\mathcal{E})\equiv
\mathcal{L}^k_a(\mathcal{E}\oplus\ldots\oplus\mathcal{E},
\mathcal{A})\equiv \mathcal{L}^k_a(\mathcal{E}, \mathcal{A}).\]It
is clear that, for every open set $U\subseteq X$, the set
$$\Omega^k_{\mathcal{A}|_U}(\mathcal{E}|_U\oplus
\ldots\oplus\mathcal{E}|_U, \mathcal{A}|_U)$$ is an
$\mathcal{A}(U)$-\textit{module}, i.e. a \textit{module over the
$\mathbb{C}$-algebra} $\mathcal{A}(U)$; hence a
$\mathbb{C}$-\textit{vector space}. Thus, based on \cite{mallios},
Proposition 1.1, p. 104 and Theorem 9.1, p. 41, we conclude that
the \textit{sheaf} $\Omega^k(\mathcal{E})$ is an
$\mathcal{A}$-\textit{module on} $X$. It follows, for every
\textit{open set} $U\subseteq X$, that
\[\Omega^k(\mathcal{E})(U)=
\Omega^k_{\mathcal{A}|_U}(\mathcal{E}|_U\oplus\ldots\oplus
\mathcal{E}|_U, \mathcal{A}|_U),\]within an
$\mathcal{A}(U)$-\textit{isomorphism} of the
$\mathcal{A}(U)$-modules concerned. In particular, one has
\[\Omega^k(\mathcal{E})(X)=
\Omega^k_\mathcal{A}(\mathcal{E}\oplus\ldots\oplus \mathcal{E},
\mathcal{A}),\]where the equality actually means an
$\mathcal{A}(X)$-\textit{isomorphism} of the
$\mathcal{A}(X)$-modules $\Omega^k(\mathcal{E})(X)$ and
$\Omega^k_\mathcal{A}(\mathcal{E}\oplus\ldots\oplus\mathcal{E},
\mathcal{A})$.

Following the classical pattern, we set \[\begin{array}{lll}
\Omega^0(\mathcal{E})= \mathcal{A}, & \mbox{and} &
\Omega^1(\mathcal{E})= \mathcal{E}^\ast,\end{array}\]as the
\textit{sheaf of germs of exterior $\mathcal{A}$-$0$-forms} and
the \textit{sheaf of germs of $\mathcal{A}$-$1$- forms}, on $X$,
respectively.

The standard \textit{exterior algebra of $k$-forms} can also be
repeated here to some significant extent. Consider for instance
the analogue, in this setting, of the usual \textit{alternation}
or \textit{anti-symmetrizer map}, \cite{ralph}, p. 101, or
\cite{greub}, p. 85, or \cite{nakahara}, p. 196, which we define
below. To this end, suppose given an $\mathcal{A}$-\textit{module
}$\mathcal{E}\equiv (\mathcal{E}, \pi, X)$ on a topological space
$X$. Let
\[\Gamma(\mathcal{E})\equiv (\Gamma(U, \mathcal{E}), \pi^U_V)\]
be the corresponding (complete) \textit{presheaf of sections} of
$\mathcal{E}$. Instead of considering the \textit{presheaf}
$T^0_k\Gamma(\mathcal{E})$, $k\in \mathbb{N}$, of
$\Gamma(\mathcal{A})$-\textit{modules} on $X$, in order to define
the \textit{alternation morphism} \[\mathbf{A}:
\mathfrak{T}^0_k(\mathcal{E})\longrightarrow
\mathfrak{T}^0_k(\mathcal{E}),\] where
$\mathfrak{T}^0_k(\mathcal{E}):=
\mathbf{S}(T^0_k\Gamma(\mathcal{E}))$, we deviate from this usual
practice to defining $\mathbf{A}$ as the $\mathcal{A}$-morphism
induced by maps
\begin{equation}\begin{array}{ll}\label{ext1}\mathbf{A}_x:
(\mathfrak{T}^0_k(\mathcal{E}))_x=
\mathcal{E}^\ast_x\otimes_{\mathcal{A}_x}\ldots\otimes_{\mathcal{A}_x}\mathcal{E}^\ast_x\longrightarrow
\mathcal{E}^\ast_x\otimes_{\mathcal{A}_x}\ldots\otimes_{\mathcal{A}_x}\mathcal{E}^\ast_x,
& x\in X,\end{array}\end{equation}such that
\[\mathbf{A}_x\mathbf{t}_x(s_{1,x},\ldots, s_{k,x})=
\frac{1}{k!}\sum_{\sigma\in S_k}(sign \sigma)\
\mathbf{t}_x(s_{1,x}, \ldots, s_{k,x}),\]where $s_{1,x}, \ldots,
s_{k,x}\in \mathcal{E}_x$, $\mathbf{t}_x:
\mathcal{E}_x\otimes_{\mathcal{A}_x}\ldots\otimes_{\mathcal{A}_x}\mathcal{E}_x\longrightarrow
\mathcal{A}_x$ is $\mathcal{A}_x$-$k$-linear, and $S_k$ is the
permutation group on $\{1, \ldots, k\}$.

The equality $(\mathfrak{T}^0_k(\mathcal{E}))_x=
\mathcal{E}^\ast_x\otimes_{\mathcal{A}_x}\ldots\otimes_{\mathcal{A}_x}\mathcal{E}^\ast_x$,
$x\in X$, holds within an $\mathcal{A}_x$-\textit{isomorphism};
for this purpose see \cite{mallios}, relation 5.9, p. 130.

The reason for this approach comes from the observation that the
$\Gamma(\mathcal{A})$-\textit{presheaf} defined by
\begin{equation}\label{exte}U\longmapsto T^0_k\Gamma(\mathcal{E})(U):=
\Gamma(U, \mathcal{E}^\ast)\otimes_{\Gamma(U,
\mathcal{A})}\ldots\otimes_{\Gamma(U, \mathcal{A})}\Gamma(U,
\mathcal{E}^\ast), \end{equation}where $U\subseteq X$ is open, along
with the restriction maps \textit{is not always complete}, cf.
\cite{mallios}, Statement 5.5, p. 129; therefore
$\mathfrak{T}^0_k(\mathcal{E})(U)$ is not always
$\mathcal{A}(U)$-isomorphic to the right hand side in the
correspondence (\ref{exte}) above, i.e., to
$T^0_k\Gamma(\mathcal{E})(U)$. Thus, in order to circumvent this
obstacle, we resort, by virtue of \cite{mallios}, Lemma 8.1, p. 36,
to anti-symmetrizers \[\begin{array}{ll} \mathbf{A}_x:
\mathfrak{T}^0_k(\mathcal{E})_x\longrightarrow
\mathfrak{T}^0_k(\mathcal{E})_x, & x\in X,\end{array}\]from which
the sought anti-symmetrizer $\mathcal{A}$-morphism $\mathbf{A}$ is
obtained.

We may now define the \textbf{exterior product} as follows.

\begin{mydf} \emph{Let $\mathcal{E}$ be a \textit{vector sheaf of rank $n$} on a
topological space $X$, and let $\xi$ and $\eta$ be elements of
$\Omega^k(\mathcal{E})$ and $\Omega^l(\mathcal{E})$, respectively.
The \textit{exterior product} of $\xi_x$ and $\eta_x$, $x\in X$,
is the \textit{germ} $\xi_x\wedge \eta_x\in
\Omega^{k+l}(\mathcal{E}_x)$, given by \[\xi_x\wedge \eta_x=
\frac{(k+l)!}{k!l!}\mathbf{A}_x(\mathcal{E}_x\otimes
\eta_x),\]that is, for all $s_{1,x},\ldots, s_{k+l, x}\in
\mathcal{E}_x$, \begin{eqnarray*}\lefteqn{\xi_x\wedge
\eta_x(s_{1,x}, \ldots, s_{k+l, x})= } \\ & &
\frac{1}{k!l!}\sum_{\sigma\in S_{k+l}}sign(\sigma)\
\xi_x(s_{\sigma(1),x},\ldots, s_{\sigma{k},x})\
\eta_x(s_{\sigma(k+1),x},\ldots,
s_{\sigma(k+l),x}).\end{eqnarray*}In particular, for $\alpha_x\in
\Omega^0(\mathcal{E}_x)\equiv \mathcal{A}_x$, $x\in X$, we put
\[\alpha_x\wedge \xi_x\equiv \xi_x\wedge \alpha_x\equiv
\alpha_x\xi_x.\]Finally, the $\mathcal{A}$-\textit{morphism}
\[\xi\wedge \eta\in \Omega^{k+l}(\mathcal{E}),\]obtained from
germs $\xi_x\wedge \eta_x$, $x\in X$, above, by virtue of
\cite{mallios}, Lemma 8.1, p. 36, is called the \textbf{exterior
product of $\xi$ and $\eta$}. Like earlier, for $\alpha\in
\Omega^0(\mathcal{E})\equiv \mathcal{A}$, we put \[\alpha\wedge
\xi\equiv \xi\wedge \alpha\equiv \alpha\xi.\]}
\end{mydf}

Note that we do not index $\wedge$, when considering the exterior
product $\xi_x\wedge \eta_x$, $x\in X$, for given $\xi\in
\Omega^k(\mathcal{E})$ and $\eta\in \Omega^l(\mathcal{E})$, in
order to avoid unnecessary meticulousness.

With this product, we define the \textbf{exterior algebra sheaf},
or the \textbf{Grassmann algebra sheaf} of the \textit{vector
sheaf} $\mathcal{E}$ of rank $n$, to be the
$\mathcal{A}$-\textit{module}
\begin{equation}\label{exte1}\Omega^\ast(\mathcal{E})\equiv \Omega^0(\mathcal{E})\oplus
\Omega^1(\mathcal{E})\oplus \ldots\oplus
\Omega^n(\mathcal{E}),\end{equation}such that
\begin{equation}\label{exte2}\Omega^\ast(\mathcal{E})_x\equiv
\Omega^0(\mathcal{E})_x\oplus \Omega^1(\mathcal{E})_x\oplus
\ldots\oplus \Omega^n(\mathcal{E})_x=
\Omega^0(\mathcal{E}_x)\oplus \Omega^1(\mathcal{E}_x)\oplus
\ldots\oplus \Omega^n(\mathcal{E}_x),\end{equation}for all $x\in
X$, where the last relation is valid, of course, within an
$\mathcal{A}_x$-\textit{isomorphism}.

The $\mathcal{A}_x$-isomorphism in (\ref{exte2}) can be obtained
in the following manner. In fact, one has
\begin{equation}\label{exte3}\Omega^1(\mathcal{E}):=
\mathcal{L}^1(\mathcal{E}, \mathcal{A})=
\mathcal{H}om_\mathcal{A}(\mathcal{E}, \mathcal{A})=
\mathcal{E},\end{equation}which implies that
\[\begin{array}{ll}\Omega^1(\mathcal{E})_x= \mathcal{E}_x, & x\in
X.\end{array}\](For the last relation in (\ref{exte3}), see
\cite{mallios}, relation (6.18), p. 136.) On the other hand, for
all $x\in X$, \[\Omega^1(\mathcal{E}_x):=
\mathcal{L}^1(\mathcal{E}_x, \mathcal{A}_x)=
\mathcal{H}om_{\mathcal{A}_x}(\mathcal{E}_x, \mathcal{A}_x)=
\mathcal{E}_x,\]because $\mathcal{E}$ is a vector sheaf of finite
rank on $X$. Thus, \[\Omega^1(\mathcal{E})_x=
\Omega^1(\mathcal{E}_x),\]for all $x\in X$. Likewise, for $k>1$,
we have \[\Omega^k(\mathcal{E})_x=:
\mathcal{L}_a^k(\mathcal{E}\oplus\ldots\oplus \mathcal{E},
\mathcal{A})_x=
\mathcal{L}^k_a(\mathcal{E}_x\oplus\ldots\oplus\mathcal{E}_x,
\mathcal{A}_x):= \Omega^k(\mathcal{E}_x).\]

Like in the classical theory, if $\alpha_i$, $i=1, \ldots, k$, are
elements of $\Omega^1(\mathcal{E})$, where $\mathcal{E}$ is a
vector sheaf of finite rank on a topological space $X$, then
\[(\alpha_{1,x}\wedge \ldots\wedge \alpha_{k,x})(s_{1,x}, \ldots,
s_{k,x})=
\sum_{\sigma}sign(\sigma)\alpha_{1,x}(s_{\sigma(1),x})\ldots\alpha_{k,x}(s_{\sigma{k},x}),\]where
$s_{i,x}:= s_i|_{\mathcal{E}_x}$, and $\alpha_{i,x}=
\alpha_i|_{\Omega^1(\mathcal{E}_x)}$, for all $X\in X$ and
$i=1,\ldots, k$, so that \[(\alpha_{1}\wedge \ldots\wedge
\alpha_{k})(s_{1}, \ldots, s_{k})=
\sum_{\sigma}sign(\sigma)\alpha_{1}(s_{\sigma(1)})\ldots\alpha_{k}(s_{\sigma{k}}).\]

\section{Skew-Symmetric $\mathcal{A}$-bilinear forms}

\begin{mydf} \emph{Let $X$ be a topological space and $(X, \mathcal{A}, \mathcal{P})$ an \textit{ordered $\mathbb{R}$-algebraized space} on $X$, (cf.\cite{mallios}, p.
316). A \textit{section} $s\in \mathcal{A}(U)$, with $U$ open in
$X$, is called \textit{strictly positive} if $s\in
\mathcal{P}(U)$, and, given $x\in U$, $s(x)\neq 0_x$.} \end{mydf}

For the purpose of what lays ahead, we need the following notion.

\begin{mydf} \emph{Let $(X, \mathcal{A})$ be an algebraized space, and
$\mathcal{E}$ an $\mathcal{A}$-module on $X$. An
$\mathcal{A}$-\textit{bilinear sheaf morphism} \[\omega\equiv
(\omega_U)_{X\supseteq U, open}: \mathcal{E}\oplus
\mathcal{E}\longrightarrow \mathcal{A}\]is called
\begin{itemize}\item \textbf{skew symmetric} provided
\[\omega_U(s, t)= -\omega_U(t, s),\]for all \textit{sections} $s, t\in
\mathcal{E}(U)$ over any \textit{open} set $U\subseteq X$. \item
$\mathbf{nondegenerate}$ if
\begin{eqnarray*}\omega_U(s, t)= 0 & \mbox{for all $t\in
\mathcal{E}(U)$, with $U$ an arbitrary open set in
$X$,}\end{eqnarray*}implies that $s= 0\in
\mathcal{E}(U)$.\end{itemize}\label{df12}}\end{mydf}

In this Definition \ref{df12}, we have identified the sheaf morphism
$\omega: \mathcal{E}\oplus \mathcal{E}\longrightarrow \mathcal{A}$
with the corresponding \textit{presheaf morphism}
$(\omega_U)_{X\supseteq U, open}: \Gamma(\mathcal{E}\oplus
\mathcal{E})\longrightarrow \Gamma(\mathcal{A})$ of
\textit{$($complete$)$ presheaves of sections}
$\Gamma(\mathcal{E}\oplus \mathcal{E})$ and $\Gamma(\mathcal{A})$.
This identification is based on the fact that, given a topological
space $X$, we have
\begin{equation}\Gamma: \mathcal{S}h_X\cong
\mathcal{C}o\mathcal{PS}h_X,\label{equivalence}\end{equation}where
$\mathcal{S}h_X$ is the category of sheaves over $X$,
$\mathcal{C}o\mathcal{PS}h_X$ is the category of complete
presheaves over $X$, and $\Gamma$ is the section functor. For
suitable details, see \cite{mallios}, Theorem 13.1, p. 73.

For the purpose of the following theorem, we assume the following
condition, referred to in the sequel as the
\textbf{inverse-positive-section condition}:
\begin{quote}{The ordered $\mathbb{R}$-algebraized $(X, \mathcal{A},
\mathcal{P})$ is such that all \textit{strictly positive sections}
of $\mathcal{A}$ are \textit{invertible}. More explicitly, if
$\mathcal{P}^*$ denotes the \textit{subsheaf} of all strictly
positive sections of $\mathcal{A}$, and $\mathcal{A}^\bullet$ the
sheaf on $X$, generated by the presheaf
\[U\longmapsto \mathcal{A}(U)^\bullet= \mathcal{A}^\bullet (U),\]
where $U$ is open in $X$, and $\mathcal{A}(U)^\bullet$ is the
\textit{group of units} of the unital $\mathbf{C}$-algebra
$\mathcal{A}(U)$, then
\begin{equation}\label{r1}\mathcal{P}^*\subset
\mathcal{A}^\bullet.\end{equation}Section-wise, (\ref{r1}) would be
understood in the following way: for any $s\in \mathcal{P}(U)$,
where $U\subseteq X$ is open, such that $s(x)\neq 0_x\in
\mathcal{A}_x$, $x\in U$, then there exists $s^{-1}\in
\mathcal{P}(U)$ such that $s\cdot s^{-1}= s^{-1}\cdot s= 1_U\in
\mathcal{A}(U)$.}\end{quote}

Furthermore, we suppose that our $\textit{ordered algebraized
space}$ $(X, \mathcal{A}, \mathcal{P})$, is also endowed with an
$\textit{absolute value}$, i.e., the following sheaf morphism,
\begin{equation}\label{absolute}|\cdot|: \mathcal{A}\longrightarrow
\mathcal{A}^+:= \mathcal{P},\end{equation}having the analogous
properties of the classical function; hence, for instance,
$\textit{the property}$ that \[|s|= \alpha\in \mathbb{R}^+\subseteq
\mathcal{A}^+(X)\Longleftrightarrow s= \pm\alpha\in
\mathbb{R}\subseteq \mathcal{A}(X).\]Now, the proof of the following
theorem is based on the classical patterns, see e.g. \cite{silva},
\cite{ralph}, \cite{berndt}, \cite{puta}, within, of course, the
present $\textit{sheaf-theoretic context}$, for which we refer to
\cite{mallios}, p. 316, Definition 8.1, along with p. 335
\textit{ff}. So, we now have the following basic result.

\begin{theo} \emph{Let $(X, \mathcal{A}, \mathcal{P}, |\cdot|)$ be an
ordered $\mathbb{R}$-algebraized space, endowed with an
\textit{absolute value}, and $\mathcal{E}$ the standard \textit{free
$\mathcal{A}$-module}, $\mathcal{A}^n$, of rank $n$ on $X$. Moreover
let $\omega: \mathcal{E}\oplus \mathcal{E}\longrightarrow
\mathcal{A}$ be a \textit{non-zero skew-symmetric} and
\textit{non-degenerate} $\mathcal{A}$-\textit{bilinear sheaf
morphism}. Then, there exists an $\mathcal{A}(U)$-\textit{basis} of
$\mathcal{A}^n(U)$, say,
\[s^U_1, \ldots, s^U_m, t^U_1, \ldots, t^U_m,\] such that
\[\begin{array}{ll} n= 2m & {} \\
\omega(s^U_i, s^U_j)= 0= \omega(t^U_i, t^U_j)& \mbox{for all $1\leq
i, j\leq m$}\\ \omega(s^U_i, t^U_j)= \delta^U_{ij} & \mbox{for all
$1\leq i, j\leq m$}.\end{array}\]} \label{the1}\end{theo}

\P With no loss of generality, we assume that $U= X$. Therefore,
since $\mathcal{A}^n\neq \{0\}$ (we already assumed that
$\mathbb{C}\equiv \mathbb{C}_X\subseteq \mathcal{A}$), there exists
an element
\[0\neq s_1\in \mathcal{A}^n(X)\cong \mathcal{A}(X)^n\](take e.g.
an element from the $\textit{canonical basis}$ of (sections) of
$\mathcal{A}^n(X)\cong \mathcal{A}(X)^n$, see \cite{mallios}, p.
123). Next, consider the ``$\textit{$\mathcal{A}(X)$-line of
$s_1$}$", i.e., \[\mathcal{A}(X)[s_1]:= \{\alpha s_1\in
\mathcal{A}^n(X):\ \alpha\in \mathcal{A}(X)\},\]which, by an obvious
$\textit{abuse of notation}$, we may still denote, for convenience,
just, by \[\mathcal{A}(s_1)\subseteq \mathcal{A}^n(X).\]Now, it is
also clear that there exists an element \[0\neq \overline{t}_1\in
\mathcal{A}^n(X)\setminus \mathcal{A}(s_1)\](just, take e.g. another
element of the previously considered basis of $\mathcal{A}^n(X)$,
$\textit{different from $s_1$}$). Furthermore, due to the hypothesis
concerning $s_1$, $\overline{t}_1$, and as well as, to that one for
$\omega$, one obtains that (see Lemma \ref{lsc})
\[\omega(s_1, \overline{t}_1)(x)\neq 0_x\in \mathcal{A}_x,\]for all $x\in X$.
Hence, based also on our
hypothesis for $\mathcal{A}$, that is, the existence of the sheaf
morphism $|\cdot|: \mathcal{A}\longrightarrow \mathcal{A}^+\equiv
\mathcal{P}$, we also obtain that,
\[|\omega(s_1, \overline{t}_1)|>0, \]that is the section
$|\omega(s_1, \overline{t}_1)|\in \mathcal{A}(X)$ is
$\textit{strictly positive}$; therefore, by assumption for $(X,
\mathcal{A})$, see the inverse-positive-section condition above,
it is also invertible in $\mathcal{A}(X)$. Hence taking further
$t_1:= u^{-1}\overline{t}_1$, with $u\equiv |\omega(s_1,
\overline{t}_1)|\in \mathcal{A}(X)$, one gets \[|\omega(s_1,
t_1)|= 1,\]which implies that \[\omega(s_1, t_1)= \pm 1\in
\mathcal{A}(X).\]Now, let us consider
\[\mathrm{S}_1:= [s_1, t_1],\]that is, the ``$\textit{flag}$"
(alias, ``$\mathcal{A}(X)$-$\textit{plane}$"), defined by $s_1$ and
$t_1$, in $\mathcal{A}^n(X)$, in effect, an
$\textit{$\mathcal{A}(X)$-module, generated by $s_1$, and $t_1$}$,
along with its $\textit{``orthogonal complement"}$ in
$\mathcal{A}^n(X)$, i.e.,
\[\mathrm{S}_1^\perp\equiv \mathrm{T}_1:= \{v\in \mathcal{A}^n(X):
\omega(v, z)= 0, \mbox{for all $z\in \mathrm{S}_1$}\}.\]Now, we
first remark that $s_1$, $t_1$ are also ``$\textit{free
generators}$" of $\mathrm{S}_1$, for, $\textit{if}$ $t_1=\alpha
s_1$, $\textit{then}$ \[\pm 1= \omega(s_1, t_1)= \omega(s_1, \alpha
s_1)= \alpha\cdot \omega(s_1, s_1)= 0,\]a $\textit{contradiction}$.
That is, $\{s_1, t_1\}$ yields  actually an
$\mathcal{A}(X)$-$\textit{basis of the flag}$ $\mathrm{S}_1$.
Furthermore, we prove that: \[\begin{array}{ll} (i) &
\mathrm{S}_1\cap \mathrm{T}_1=\{0\}, \ \mbox{and}\\ (ii) &
\mathrm{S}_1+ \mathrm{T}_1= \mathcal{A}^n(X).\end{array}\]Indeed,
$(i)$ if $z\equiv \alpha s_1+ \beta t_1\in \mathrm{S}_1\cap
\mathrm{T}_1$, with $\alpha, \beta\in \mathcal{A}(X)$, one gets, by
the very definition of $\mathrm{S}_1$, $\mathrm{T}_1$, and the fact
that $\omega(s_1, t_1)= 1$, that, \[\begin{array}{ll} \omega(z,
s_1)= \beta= 0, & \mbox{and $\omega(z, t_1)=\alpha=
0$},\end{array}\]that is, $z=0$, $\textit{which proves}$ $(i)$. On
the other hand, $(ii)$ $\textit{for every}$ $z\in \mathcal{A}^n(X)$,
one has, \[ z= (-\omega(z, s_{1})t_{1}+ \omega(z, t_{1})s_{1})+ (z+
\omega(z, s_{1})t_{1}- \omega(z, t_{1})s_{1}),\]with
\[-\omega(z, s_{1})t_{1}+ \omega(z, t_{1})s_{1}\in
\mathrm{S}_{1},\] and \[z+ \omega(z, s_{1})t_{1}-\omega(z,
t_{1})s_{1}\in \mathrm{T}_{1}.\] Thus, \[\mathcal{A}^n(X)=
\mathrm{S}_{1}\oplus \mathrm{T}_{1}.\]

Now, in a manner similar to the manner we found the elements
$s_1$, $t_1\in \mathrm{S}_1$ with $\omega(s_1, t_1)= 1$, we
conclude the existence of $\textit{sections }$ $s_2, t_2\in
\mathrm{T}_1\setminus \{0\}$, $\textit{such that }$ \[\omega(s_2,
t_2)= 1\in \mathcal{A}(X);\]while we further consider the flag
\[\mathrm{S}_2:= [s_2, t_2],\]along with \[\mathrm{T}_2\equiv
\mathrm{S}_2^\perp:= \{v\in \mathcal{A}^n(X):\ \omega(v, z)= 0,\
z\in \mathrm{S}_2\}.\]Yet, we still prove in a similar way, as
before, that \[\mathrm{T}_1= \mathrm{S}_2\oplus \mathrm{T}_2,\]so
that one obtains, \[\mathcal{A}^n(X)= \mathrm{S}_1\oplus
\mathrm{S}_2\oplus \mathrm{T}_2,\]and so on. Now, the above
process stops eventually, due to the $\textit{finite rank of}$
$\mathcal{A}^n(X)$, so that one finally obtains
\[\mathcal{A}^n(X)= \mathrm{S}_{1}\oplus \mathrm{S}_{2}\oplus
\ldots\oplus \mathrm{S}_{m}\]with the generators, $s_{i}$,
$t_{i}$, of $\mathrm{S}_{i}$ ($1\leq i\leq n$) having the property
that
\[\begin{array}{ll}
\omega(s_{i}, s_{j})=  0  =  \omega(t_{i}, t_{j}) & \\
\omega(s_{i}, t_{j}) = \delta_{ij}. &
\end{array}\] Hence, the proof is
finished. \proofterminator

In the proof of the previous Theorem \ref{the1}, one still
essentially applies the following standard fact of the classical
theory, which for convenience we also formulate, within the present
context:
\begin{lem}\emph{Let $\mathcal{E}$ be a free $\mathcal{A}$-\textit{module} of rank
$n\in \mathbb{N}$ on a topological space $X$. Then, a family
$\{s_i\}_{i\in I}$ of \textit{global sections} of $\mathcal{E}$,
i.e., $\{s_i\}_{i\in I}\subseteq \mathcal{E}(X)$, is
$\mathcal{A}(X)$-\textit{linearly independent} if, and only if, the
relation
\[\sum_{i\in I}\alpha_is_i= 0,\]with $\{\alpha_i\}_{i\in I}\subseteq
\mathcal{A}(X)$, having finite support, implies $\alpha_i= 0$, for
any $i\in I$. \label{lsc}}\end{lem} For the proof of Lemma
\ref{lsc}, one can follow, for instance, the analogous argument in
\cite{fbourbaki}, Chap II; p. 25, remarks after Definition 10. Yet,
for convenience, we recall that the term ``\textit{finite support}"
for the family $\{\alpha_i\}_{i\in I}\subseteq \mathcal{A}(X)$,
means that $\alpha_i\neq 0$, \textit{only for finitely many indices}
$i\in I$, \textit{and the rest of the $\alpha_i$ being} $0$; so that
the \textit{sum} used above acquires then a meaning.

When the skew-symmetric $\mathcal{A}$-bilinear sheaf morphism
$\omega$ is not necessarily nondegenerate, then in place of Theorem
\ref{the1}, we have the theorem below. Let us first give the
following definition:

\begin{mydf} \emph{Let $\omega\equiv (\omega^U): \mathcal{A}^n\oplus
\mathcal{A}^n\longrightarrow \mathcal{A}$ be an
$\mathcal{A}$-bilinear morphism on the standard free
$\mathcal{A}$-module $\mathcal{A}^n$. The \textit{rank} of $\omega$
is the rank of the matrix $(\omega_{ij}^U)$, with $\omega_{ij}^U=
\omega^U(\varepsilon_i^U, \varepsilon_j^U)$, where
$\{\varepsilon_i^U\}_{1\leq i\leq n}$ is the \textit{Kronecker
gauge} of $\mathcal{A}^n(U)$, and $U$ is any \textit{open subset} of
$X$.}
\end{mydf} By the methods of Linear
Algebra, one can easily show that the rank of $\omega$ is
independent of the basis considered.

Note the notation $\omega\equiv (\omega^U)$ instead of the usual
$\omega\equiv (\omega_U)$; the reason being that we want to avoid,
in the sequel, too many sub-indices.

\begin{theo} \emph{Let $(X, \mathcal{A}, \mathcal{P}, |\cdot|)$ and $\mathcal{E}$ be as
in Theorem \ref{the1}. Let $\omega: \mathcal{E}\oplus
\mathcal{E}\longrightarrow \mathcal{A}$ be a skew-symmetric
$\mathcal{A}$-bilinear sheaf morphism of rank $r$. Then, $r=2m$, for
some integer $m$, and for every $x\in X$ there are an open
neighborhood $U\subseteq X$ of $x$ and a basis \[s_1^U,\ldots,
s_n^U\in \mathcal{A}^n(U)= \mathcal{A}(U)^n\]such that the matrix of
$\omega^U\equiv \omega_U$ is \[\left[\begin{array}{ccc} 0 &
\mbox{I}_m & 0\\ -\mbox{I}_m & 0 & 0\\ 0 & 0 &
0\end{array}\right],\]where $\omega_U$ is the component of $\omega$
relative to $U$.} \end{theo} \P Fix $x\in X$. Because of the
continuity of $\omega_y$, $y\in X$, and of the fact that the rank of
$\omega$ is $r$, there exist an open neighborhood $U_1$ of $x$ and
$s_1^{U_1}, \bar{s}_{m+1}^{U_1}\in \mathcal{A}^n(U_1)$ such that
\[\omega^{U_1}(s_1^{U_1}, \bar{s}_{m+1}^{U_1})(y)\neq 0_y\in
\mathcal{A}_y,\]for all $y\in U_1$. Now, put $t_{1,m+1}^{U_1}:=
\omega^{U_1}(s_1^{U_1}, \bar{s}_{m+1}^{U_1})\in \mathcal{A}(U_1)$.
Based on our hypothesis for $\mathcal{A}$, we have that $u_{1,m+1}:=
|t_{1,m+1}^{U_1}|>0$; therefore by the inverse-positive-section
condition, $u_{1,m+1}\in \mathcal{A}^\bullet(U_1)=
\mathcal{A}(U_1)^\bullet$. Hence, taking further $s^{U_1}_{m+1}=
u^{-1}_{1,m+1}\bar{s}^{U_1}_{m+1}$, one gets
\[|\omega^{U_1}(s_1^{U_1}, s^{U_1}_{m+1})|= 1.\]Assuming that $\omega^{U_1}(s_1^{U_1}, s^{U_1}_{m+1})|=
1_{U_1}=:1$, the matrix of $\omega^{U_1}$ in the
$\mathcal{A}(U_1)$-module $\mathrm{S}_1^{U_1}:= [s_1^{U_1},
s_{m+1}^{U_1}]$, that is, the $\mathcal{A}(U_1)$-module spanned by
$s_1^{U_1}$ and $s_{m+1}^{U_1}$, is \[\left[\begin{array}{cc}0 & 1\\
-1 & 0\end{array}\right].\]

Let $(\mathrm{S}_1^{U_1})^\bot$ be the $\omega^{U_1}$-orthogonal
complement of $\mathrm{S}_1^{U_1}$ in $\mathcal{A}^n(U_1)$. As in
the proof of Theorem \ref{the1}, one shows that
\[\mathcal{A}^n(U_1)= \mathrm{S}_1^{U_1}\oplus
(\mathrm{S}_1^{U_1})^\bot.\]Now, we repeat the process on
$(\mathrm{S}_1^{U_1})^\bot$ to get an open neighborhood
$U_2\subseteq X$ of $x$, and $s_2^{U_2}$ and $s_{m+2}^{U_2}$ such
that $\omega^{U_2}(s_2^{U_2}, s_{m+2}^{U_2})= 1$ and continue
inductively. Eventually, this process will stop because the rank of
$\omega$ is finite. Certainly by taking $U= \cap_{i=1}^mU_i$, where
$r= 2m$, one sees that $\omega^U$ has the stated matrix in the basis
$\{s_1^U, \ldots, s_n^U\}$. \proofterminator

Note that if we denote by $\{(s_i^U)^\ast\}_{1\leq i\leq m}$ the
dual basis of $\{s_i^U\}_{1\leq i\leq m}$, then
\begin{equation}\omega^U= \sum_{i=1}^m(s_i^U)^\ast\wedge
(s_{m+i}^U)^\ast.\label{nondeg}\end{equation}

\begin{cor} \emph{Let $(X, \mathcal{A}, \mathcal{P})$ be an ordered $\mathbb{R}$-algebraized space such that $\mathcal{P}^*\subseteq
\mathcal{A}^\bullet$. Let $\mathcal{E}$ be a \textit{vector sheaf of
rank} $n$ on $X$, and $\omega: \mathcal{E}\oplus
\mathcal{E}\longrightarrow \mathcal{A}$ a \textit{skew-symmetric}
and \textit{non-degenerate} $\mathcal{A}$-\textit{bilinear
morphism}. Then, given an \textit{open} neighborhood $U$ of $x\in X$
such that $\mathcal{E}|_U= \mathcal{A}^n|_U$, where $x\in X$ is
arbitrary, there exists a \textit{basis}, $s_1^U, \ldots, s_m^U,
t_1^U,\ldots, t_m^U$, of $\mathcal{A}^n(U)$ such that
\[\begin{array}{ll} n=2m & {}\\ \omega^U(s_i^U, s_j^U)=0 =
\omega^U(t_i^U, t_j^U) & \mbox{for all $1\leq i, j\leq m$}\\
\omega^U(s_i^U, t_j^U)= \delta_{ij}^U & \mbox{for all $1\leq i,
j\leq m$}.\end{array}\]The pair $(\mathcal{E}, \omega)$ is called a
\textbf{locally free symplectic $\mathcal{A}$-module}, alias
\textbf{symplectic vector sheaf}, of rank $n$, on $X$.}\end{cor} \P
One proceeds in the same manner as for the proof of Theorem
\ref{the1}, the small nuance being that one works locally.
\proofterminator

\begin{mydf} \emph{Let $(X, \mathcal{A}, \mathcal{P})$ be an ordered $\mathbb{R}$-algebraized space, satisfying the
inverse-positive-section condition. The non-degenerate
skew-symmetric $\mathcal{A}$-bilinear morphism $\omega:
\mathcal{A}^n\oplus \mathcal{A}^n\longrightarrow \mathcal{A}$ is
called a \textbf{linear symplectic $\mathcal{A}$-structure} on the
standard free $\mathcal{A}$-module $\mathcal{A}^n$, and the pair
$(\mathcal{A}^n, \omega)$ is called a \textbf{$($standard$)$ free
symplectic $\mathcal{A}$-module}. More generally, let $\mathcal{E}$
be an $\mathcal{A}$-module on $X$, and $\omega: \mathcal{E}\times
\mathcal{E}\longrightarrow \mathcal{A}$ a non-degenerate,
skew-symmetric $\mathcal{A}$-bilinear morphism. Then, the pair
$(\mathcal{E}, \omega)$ is called a \textbf{symplectic
$\mathcal{A}$-module} on $X$.}
\end{mydf}

With respect to the notation of Theorem \ref{the1}, the basis,
$s_1,\ldots, s_m, t_1, \ldots, t_m$, is called a $\textbf{symplectic
basis}$ of the standard symplectic free $\mathcal{A}$-module
$(\mathcal{A}^n, \omega)$. It is clear that, with respect to a
symplectic basis, the matrix representing $\omega$ is, as in the
classical case, given by \begin{equation}J= \left[\begin{array}{ll}
0 & \mbox{I}_n \\ \hspace{-2mm}-\mbox{I}_n &
0\end{array}\right].\label{matrix}\end{equation}

\begin{exam} \emph{Let $\mathcal{E}$ be a vector sheaf of rank $n$ on $X$,
and consider the direct sum $\mathcal{E}\oplus \mathcal{E}^\ast$.
The $\mathcal{A}$-bilinear morphism \[\omega: (\mathcal{E}\oplus
\mathcal{E}^\ast)\oplus (\mathcal{E}\oplus
\mathcal{E}^\ast)\longrightarrow \mathcal{A},\]defined by
\[\omega^U\left((s_1^U, \alpha_1^U), (s_2^U, \alpha_2^U)\right)=
\alpha_2^U(s_1^U)- \alpha_1^U(s_2^U),\]where $U$ is a local gauge of
$\mathcal{E}$, and $s_i^U\in \mathcal{E}(U)= \mathcal{A}^n(U)=
\mathcal{A}(U)^n$ and $\alpha_i^U\in \mathcal{E}^\ast(U)=
(\mathcal{A}^n)(U)= \mathcal{A}^n(U)$, is a symplectic
$\mathcal{A}$-morphism. This example shows that any vector sheaf
$\mathcal{F}$ of \textit{even} rank admits a symplectic
$\mathcal{A}$-structure, for one has \[\mathcal{F}|_U=
\mathcal{A}^{2n}|_U= \mathcal{A}^n|_U\oplus \mathcal{A}^n|_U=
\mathcal{A}^n|_U\oplus (\mathcal{A}^n)^\ast|_U,\]where $U$ is a
local gauge of $\mathcal{F}$. For the last equality of the previous
line, see [\cite{mallios}, relation (3.14), p.122].}
\end{exam}

We now would like to show that a useful criterion for non-degeneracy
is also available in this setting. To this effect, we suppose as
usual that our $\textit{ordered alge-}$ $\textit{braized space}$
$(X, \mathcal{A}, \mathcal{P})$ is enriched with absolute value, see
(\ref{absolute}), and $\textit{square root}$; the latter means a
$\textit{morphism of (complete) presheaves}$
\[\begin{array}{lll}\sqrt{}: \Gamma(\mathcal{P})\longrightarrow
\Gamma(\mathcal{P}), &  s\longmapsto \sqrt{s}, & \mbox{for all $s\in
\mathcal{P}(U)$},\end{array}\] with $U$ running over all the open
subsets of $X$. In addition, we assume that the pair $(\mathcal{A},
\rho)$ is a \textit{Riemann} $\mathcal{A}$-module on $X$, see
\cite{mallios}, p. 320. We assume as well that
$\overline{\rho}\equiv \rho^n$ (cf. \cite{mallios}, p. 324) is the
extension of $\rho$ to the standard free $\mathcal{A}$-module
$\mathcal{A}^n$ on $X$.

\begin{mydf} \emph{Let $(X, \mathcal{A}, \mathcal{P})$ and
$\overline{\rho}$ be as above, and $\{s_i\}_{1\leq i\leq n}\subseteq
\mathcal{A}^n(X)= \mathcal{A}(X)^n$ be a basis of
$\mathcal{A}^n(X)$. Denoting by
\[\bigwedge^n(\mathcal{A}^n)^\ast\]the $n$-th \textit{exterior power} of
the $\mathcal{A}$-module $(\mathcal{A}^n)^\ast$, see
$\cite{mallios}$, p. 307, the \textit{section}
\begin{equation}\label{vol}\Omega =\sqrt{|\det(\overline{\rho}(s_i,
s_j))|} s_1^\ast\wedge \ldots\wedge s_n^\ast\in
\left(\Lambda^n(\mathcal{A}^n)^\ast\right)(X),\end{equation} where
$\{s_i^\ast\}_{1\leq 1\leq n}\subseteq (\mathcal{A}^n)^\ast(X)$ is
the dual $\mathcal{A}$-basis of $\{s_i\}_{1\leq i\leq n}$, is called
a $\textbf{volume element of the $\mathcal{A}$-metric
$\overline{\rho}$}$. In the sequel, for the sake of brevity, the
scaling factor $\sqrt{|\det(\overline{\rho}(s_1, s_j))|}$, above,
will be denoted by $\sqrt{|\overline{\rho}(s)|}\equiv
\sqrt{|\overline{\rho}|}$, where $s\equiv \{s_1,\ldots, s_n\}$.}
\label{nvol}\end{mydf}

It is clear that for an orthonormal gauge $\{s_i\}_{1\leq i\leq
n}$ of $\mathcal{A}^n(X)$ (see \cite{mallios}, p. 340), relation
(\ref{vol}) becomes  \[\Omega= s_1^\ast\wedge \ldots\wedge
s_n^\ast.\]

Like in the classical case, we have:

\begin{cor} \emph{Let $\omega: \mathcal{A}^n\oplus
\mathcal{A}^n\longrightarrow \mathcal{A}$ be a
$\mathcal{A}$-bilinear and skew-symmetric $\mathcal{A}$-morphism on
the standard free $\mathcal{A}$-module $\mathcal{A}^n$, where $(X,
\mathcal{A})$ is an enriched $($with square root and absolute
value$)$ ordered algebraized space and $(\mathcal{A}, \rho)$ a
Riemannian $\mathcal{A}$-module. Then, $\omega$ is
\textit{non-degenerate} if and only if $n=2m$, for some $m\in
\mathbb{N}$, and $\omega^m= \omega\wedge \ldots\wedge \omega$ is a
\textit{volume element} on $\mathcal{A}^n$.}
\end{cor} Note that \[\omega^m\equiv ({\omega^U}^m):=
({\omega^U}\wedge\ldots\wedge {\omega^U}).\] Therefore, that
$\omega^m$ is a volume element on $\mathcal{A}^n$ means that
${\omega^U}\wedge \ldots\wedge {\omega^U}$ is a volume element on
$\mathcal{A}^n(U)$, for every open $U\subseteq X$.

\P We refer to the proof in Abraham-Marsden [\cite{ralph}, p. 165]
for detail.

Suppose that $\omega$ is non-degenrate. By Theorem \ref{the1},
$n=2m$, for some $m\in \mathbb{N}$. By virtue of Equation
(\ref{nondeg}) and of induction, one sees that \[{\omega^U}^m= m!
(-1)^{[m/2]}{s_1^U}^\ast\wedge \ldots\wedge {s_{2m}^U}^\ast,\]where
$[m/2]$ is the largest integer in $m/2$. Thus, assuming that
$\sqrt{|\overline{\rho}|}= m! (-1)^{[m/2]}$, ${\omega^U}^m$ is a
volume. The converse is clear. \proofterminator

Following Abraham-Marsden [\cite{ralph}, p. 167], we set that
given a free symplectic $\mathcal{A}$-module $(\mathcal{A}^n,
\omega)$, the volume element \[\Omega_\omega=
\frac{(-1)^{[m/2]}}{m!}\omega^m\] defines an
$\textbf{orientation}$ on $\mathcal{A}^n$.

In all that precedes, by the $\mathcal{A}$-bilinear morphism
$\omega: \mathcal{E}\oplus \mathcal{E}\longrightarrow
\mathcal{A}$, we meant the map on sections of the corresponding
$\mathcal{A}$-modules. So, using the presheaf $(\Gamma(U,
\mathcal{A}^n), \sigma^U_V)$ of sections of the free
$\mathcal{A}$-module $\mathcal{A}^n$ of Theorem \ref{the1}, with
$U$ ranging over the topology $\mathcal{T}$ of $\mathcal{A}^n$, we
adopt the following classical terminology. Let $S$ be an
$A(X)$-submodule of the $A(X)$-module $\Gamma(X, \mathcal{A}^n)$.
Then,
\begin{itemize} \item $S$ is called \textit{symplectic} if $\omega|_S$ is
non-degenerate. For example, $\mathrm{S}_1$, in the proof of Theorem
\ref{the1}, is symplectic. \item $S$ is called \textit{isotropic} if
$\omega|_S\equiv 0$. For instance, the span of $s_1$, $s_2$, in
Theorem \ref{the1} is isotropic.\end{itemize}

\begin{mydf} \emph{Let $(\mathcal{E}, \omega)$ and $(\mathcal{E}',
\omega')$ be symplectic $\mathcal{A}$-modules on the same
topological space $X$. An $\mathcal{A}$-morphism $\varphi:
\mathcal{E}\longrightarrow \mathcal{E}'$ is called
$\textbf{symplectic}$ if
\begin{equation}\varphi^*\omega':=\omega'\circ (\varphi\times
\varphi)= \omega \label{symplecto1},\end{equation} that is, for any
$s, t\in \mathcal{E}(U)$, $U\subseteq X$ open,
\begin{equation}(\varphi^\ast\omega')(s, t):= \omega'(\varphi(s),
\varphi(t))= \omega(s, t).\label{symplecto}\end{equation} A
symplectic $\mathcal{A}$-isomorphism is called an $
\textbf{$\mathcal{A}$-symplectomorphism}$. Symplectic
$\mathcal{A}$-modules $(\mathcal{E}, \omega)$ and $(\mathcal{E}',
\omega')$ are called $\mathcal{A}$-symplectomorphic if there is an
$\mathcal{A}$-symplectomorphism $\varphi$ between them.}
\label{symplectomorphism}\end{mydf}

Strictly speaking, Equations (\ref{symplecto}) and
(\ref{symplecto1}) should be written as follows:
\[(\varphi^\ast_U\omega'_U)(s, t):= \omega'_U(\varphi_U(s),
\varphi_U(t))= \omega_U(s, t), \] and
\[\varphi_U^\ast(\omega'_U):= \omega'_U\circ (\varphi_U\times \varphi_U)=
\omega_U,\] respectively, where $U$ varies over the topology of
$X$.

It is clear that if $(\mathcal{E}, \omega)$ and $(\mathcal{E}',
\omega')$ are symplectomorphic, then they are of the same rank,
and their rank is an even positive integer. It is also clear that
the set of symplectic $\mathcal{A}$-modules, defined over the same
topological space, can be partitioned into equivalence classes.
Furthermore, \textit{a symplectic $\mathcal{A}$-morphism is}
necessarily \textit{injective}, since if $\varphi:=
(\varphi_V)_{X\supseteq V, open}$ and $\varphi_V(s)= 0$, where
$s\in \Gamma(V, \mathcal{E})$, then necessarily $s=0$, for
$\omega$ is non-degenerate.

\begin{lem} \emph{Let $\mathcal{A}$ be a unital $\mathbb{C}$-algebra
sheaf on a topological space $X$, and let
\[\mathcal{S}p(\mathcal{E}, \omega)\equiv
\mathcal{S}p\ \mathcal{E}\]be the sheaf on $X$, generated by the
presheaf \begin{equation}U\longmapsto (\mathrm{S}p\
\mathcal{E})(U),\label{symplecto2}\end{equation}where $U$ varies
over the topology of $X$, such that for every open set $U\subseteq
X$, $(\mathrm{S}p\ \mathcal{E})(U)$ is the group $($under
composition$)$ of all $\mathcal{A}|_U$-symplectomorphisms
\[(\mathcal{E}|_U, \omega|_U\equiv \omega)\longrightarrow
(\mathcal{E}|_U, \omega|_U\equiv \omega).\]Then, the
correspondence, given by $(\ref{symplecto2})$, yields a complete
presheaf of groups on $X$; so that one obtains
\[(\mathcal{S}p\ \mathcal{E})(U)= (\mathrm{S}p\ \mathcal{E})(U),\]up
to a group isomorphism, for every open set $U\subseteq X$. The sheaf
$\mathcal{S}p\ \mathcal{E}$ is called the \textbf{symplectic group
sheaf}, or even \textbf{group sheaf of symplectomorphisms} of
$\mathcal{E}$ $($ in fact, of $(\mathcal{E}, \omega))$ on $X$.}
\label{lem21}\end{lem}

\P We first show that for all open set $U\subseteq X$,
$(\mathrm{S}p\ \mathcal{E})(U)$ is a group. In fact, let us
consider the subset
\[\mathrm{G}\mathrm{L}_{\mathcal{A}|_U}(\mathcal{E}|_U,
\mathcal{E}|_U)\subseteq
\mathrm{H}om_{\mathcal{A}|_U}(\mathcal{E}|_U, \mathcal{E}|_U)\]of
all invertible elements of the $\mathcal{A}(U)$-module
$\mathrm{H}om_{\mathcal{A}|_U}(\mathcal{E}|_U, \mathcal{E}|_U)$.
For $\mathcal{A}|_U$-morphisms $\varphi$, $\psi\in
\mathrm{G}\mathrm{L}_{\mathcal{A}|_U}(\mathcal{E}|_U,
\mathcal{E}|_U)$, we define \[\varphi\circ \psi= (\varphi_V\circ
\psi_V)_{U\supseteq V, open},\]and \[\varphi^{-1}=
(\varphi^{-1}_V)_{U\supseteq V, open}.\]It is easy to see that
under the above law of composition
$\mathrm{G}\mathrm{L}_{\mathcal{A}|_U}(\mathcal{E}|_U,
\mathcal{E}|_U)$ forms a group. Therefore, in order to show that
$(\mathrm{S}p\ \mathcal{E})(U)$, where $U$ is an open subset of
$X$, is a group, we need only show that if $\varphi, \psi\in
(\mathrm{S}p\ \mathcal{E})(U)$, then $\varphi\circ \psi$ and
$\varphi^{-1}\in (\mathrm{S}p\ \mathcal{E})(U)$. To this end, let
$V$ be an open subset of $U$. Then, we have
\[(\varphi_V\circ \psi_V)^\ast(\omega)=
(\psi^\ast_V\circ \varphi^\ast_V)(\omega)= \psi^\ast_V(
\varphi^\ast_V(\omega)) = \psi^\ast_V(\omega)= \omega ,\]from
which we deduce that $\varphi\circ \psi\in (\mathrm{S}p\
\mathcal{E})(U)$.

On the other hand, one also obtains
\[(\varphi^{-1}_V)^\ast(\omega) =
(\varphi^{-1}_V)^\ast(\varphi_V^\ast(\omega)) = (\varphi_V\circ
\varphi^{-1}_V)^\ast(\omega)= \omega,\]so that $\varphi^{-1}\in
(\mathrm{S}p\ \mathcal{E})(U)$, as well.

Let us now show that $(\ref{symplecto2})$ yields a complete
presheaf of groups. It is easy to see that Correspondence
$(\ref{symplecto2})$, along with the obvious restriction maps,
defines a presheaf of groups on $X$. Thus, we just prove that the
presheaf of groups defined, on $X$, by $(\ref{symplecto2})$ is
complete.

Indeed, let $U$ be an open subset of $X$ and $\mathcal{U}=
\{U_\alpha\}_{\alpha\in I}$ an open covering of $U$; let $\varphi,
\psi\in (\mathrm{S}p\ \mathcal{E})(U)$ such that
\[\rho^U_{U_\alpha}(\varphi)\equiv \varphi_{U_\alpha}:=
\varphi_\alpha= \psi_\alpha=: \psi_{U_\alpha}\equiv
\rho^U_{U_\alpha},\] for all $\alpha\in I$, and where the
$\rho^U_{U_\alpha}$ are the restriction maps characterizing the
presheaf $((\mathrm{S}p\ \mathcal{E})(U), V)$. Since
\[\begin{array}{ll}(\mathrm{S}p\ \mathcal{E})(U)\subseteq
\mathrm{G}\mathrm{L}_{\mathcal{A}|_U}(\mathcal{E}|_U,
\mathcal{E}|_U), &
\mathrm{G}\mathrm{L}_{\mathcal{A}|_U}(\mathcal{E}|_U,
\mathcal{E}|_U)\subseteq
\mathrm{H}om_{\mathcal{A}|_U}(\mathcal{E}|_U,
\mathcal{E}|_U),\end{array}\] and the
$\{\rho^U_{U_\alpha}\}_{\alpha\in I}$ are also the restriction
maps making the diagram
\begin{equation}\label{symplecto3} U\longmapsto
\mathrm{H}om_{\mathcal{A}|_U}(\mathcal{E}|_U,
\mathcal{E}|_U)\end{equation} into a presheaf, it follows that
$\varphi= \psi$. So the presheaf on $X$, given by
$(\ref{symplecto2})$, satisfies Condition $(S1)$ of presheaves,
see \cite{mallios}, p. 46.

For axiom $(S2)$, see \cite{mallios}, p. 47, let
\[(\varphi_\alpha)_{\alpha\in I}\in \prod_{\alpha\in
I}(\mathrm{S}p\ \mathcal{E})(U_\alpha)\subseteq \prod_{\alpha\in
I}\mathrm{H}om_{\mathcal{A}|_{U_\alpha}}(\mathcal{E}|_{U_\alpha},
\mathcal{E}|_{U_\alpha})\] be such that
\[\rho^{U_\alpha}_{U_\alpha\cap U_\beta}(\varphi_\alpha)\equiv
\varphi_\alpha|_{U_\alpha\cap U_\beta}=
\varphi_\beta|_{U_\alpha\cap U_\beta}\equiv
\rho^{U_\beta}_{U_\alpha\cap U_\beta}(\varphi_\beta)\]for any
$\alpha, \beta\in I$, with $U_\alpha\cap U_\beta\neq \emptyset$.
Hence, since $(\ref{symplecto3})$ yields a complete presheaf,
there exists an element $\varphi\in
\mathrm{H}om_{\mathcal{A}|_U}(\mathcal{E}|_U, \mathcal{E}|_U)$
such that one has \[\begin{array}{ll}\varphi|_{U_\alpha}=
\varphi_\alpha, & \alpha\in I.\end{array}\]It only remains to show
that \[\varphi^\ast\omega= \omega,\]where $\omega\equiv \omega|_U$
is a symplectic structure on $\mathcal{E}|_U$.

To this end, we first observe that $\varphi^\ast\omega$,
$\omega\in \mathrm{H}om_{\mathcal{A}|_U}((\mathcal{E}\times
\mathcal{E})|_U, \mathcal{A}|_U)$, with \[U\longmapsto
\mathrm{H}om_{\mathcal{A}|_U}((\mathcal{E}\times \mathcal{E})|_U,
\mathcal{A}|_U)\] defining a complete presheaf of
$\mathcal{A}$-modules on $X$, see \cite{mallios}, p. 134. But,
\[\varphi^\ast\omega|_{U_\alpha}=
(\varphi|_{U_\alpha})^\ast\omega|_{U_\alpha}=
\varphi^\ast_\alpha\omega|_{U_\alpha}= \omega|_{U_\alpha},\]
therefore
\[\varphi^\ast\omega= \omega,\] as desired. \proofterminator

On the other hand, the preceding notion of $\textit{symplectic
sheaf of groups}$ can also be defined $\textit{through}$ an
application of the $\textit{isomorphism}$ $\Gamma$, which is given
in $(\ref{equivalence})$. Namely, since by $(\ref{equivalence})$,
$\varphi\equiv (\varphi_U): \Gamma(\mathcal{E})\longrightarrow
\Gamma(\mathcal{E})$ is an $A$-symplectomorphism $\textit{if and
only if}$ the corresponding $\mathcal{A}$-morphism is
symplectomorphic, the symplectic sheaf of groups can be viewed as
consisting of all $\mathcal{A}$-$\textit{symplectomorphisms}$
\[\varphi: \mathcal{E}\longrightarrow \mathcal{E}.\]

Now, suppose that $\mathcal{E}$ is the $\textit{standard free
$\mathcal{A}$-module}$ $\mathcal{A}^{2n}$ of rank $2n$, and let
$\omega: \mathcal{A}^{2n}\oplus \mathcal{A}^{2n}\longrightarrow
\mathcal{A}$ be a $\textit{skew-symmetric, non-degenerate}$
$\mathcal{A}$-$\textit{bilinear morphism}$. Let $\{s_i\}_{1\leq
i\leq 2n}$ be a $\textit{basis}$ of $\mathcal{A}^{2n}(X)$ such
that $\textit{Theorem \ref{the1} holds}$, and let $\varphi\in
\mathcal{S}p\ \mathcal{A}^{2n}(X)$. Let's consider the
$\textit{full matrix algebra sheaf}$ $M_{2n}$ (see \cite{mallios},
p. 280) induced by the presheaf \[U\longmapsto
M_{2n}(\mathcal{A}(U)),\]where $U\subseteq X$ is $\textit{open}$,
and the range, $M_{2n}(\mathcal{A}(U))$, of the preceding map
consists of all $2n\times 2n$-matrices with entries in the
$\mathbb{C}$-\textit{algebra (unital and commutative)}
\[\mathcal{A}(U)\equiv \Gamma(U, \mathcal{A}).\]
Since
\[\mathcal{S}p\ \mathcal{A}^{2n}(X)\subseteq
\mathcal{H}om_\mathcal{A}(\mathcal{A}^{2n},
\mathcal{A}^{2n})(X),\]and by Statement 3.17, \cite{mallios}, p.
293, \[\mathcal{H}om_\mathcal{A}(\mathcal{A}^{2n},
\mathcal{A}^{2n})(X)= M_{2n}(\mathcal{A}(X))=
M_{2n}(\mathcal{A})(X),\]it follows that if $M$ is the
$\textit{matrix representing}$ $\varphi\in \mathcal{S}p\
\mathcal{A}^{2n}(X)$ with respect to the basis $\{s_i\}_{1\leq
i\leq 2n}$ above, then Equation (\ref{symplecto1}) becomes in
matrix form \begin{equation} {}^tMJM=
J,\label{unit}\end{equation}where $J$ is the matrix
$(\ref{matrix})$.

From (\ref{unit}), we deduce the following corollary.

\begin{cor} \emph{The \textit{determinant} of an
$\mathcal{A}$-\textit{symplectomorphism} $$\varphi:
(\mathcal{A}^{2n}, \omega)\longrightarrow (\mathcal{A}^{2n},
\omega),$$ where $\mathcal{A}$ is a \textit{unital and
commutative} $\mathbb{C}$-\textit{algebra sheaf} on a topological
space $X$,is the \textit{global identity section} $1\in \Gamma(X,
\mathcal{A}^{2n})$. More explicitly, given $M$ as the matrix
representing the $\mathcal{A}$-symplectomorphism $\varphi$, as in
the paragraph preceding the corollary, one has
\begin{equation}\label{cor24}\overline{{\partial}et}(M)= 1\in \Gamma(X,
\mathcal{A}^{2n}).\end{equation}\medskip We refer to
$\cite{mallios}$, p. 294, for the definition of the determinant
morphisms $\partial et$ and $\overline{{\partial}et}$.}
\end{cor}

\section{The Characteristic Polynomial Section}Let $\mathcal{A}$ and
$\mathcal{B}$ be \textit{sheaves of algebras} on a topological space
$(X, \mathcal{T})$, and let $\varphi: \mathcal{A}\longrightarrow
\mathcal{B}$ be a sheaf morphism such that, if
\[\overline{\varphi}\equiv (\overline{\varphi}_U)_{U\in
\mathcal{T}}: \Gamma(\mathcal{A})\longrightarrow
\Gamma(\mathcal{B})\] is the corresponding morphism between the
associated complete presheaves of sections $\Gamma(\mathcal{A})$ and
$\Gamma(\mathcal{B})$, then, for all $U\in \mathcal{T}$,
\begin{equation}\label{center}
\overline{\varphi}_U(\mathcal{A}(U))\subseteq
C(\mathcal{B}(U)),\end{equation}where $C(\mathcal{B}(U))$ stands for
the center of the ring $\mathcal{B}(U)$, cf. [\cite{lang}, p.121] .
Explicitly, Equation (\ref{center}) means that
\[\overline{\varphi}_U(s)t= t\overline{\varphi}_U(s)\]for all $s\in
\mathcal{A}(U)$ and $t\in \mathcal{B}(U)$. Now, given $s\in
\mathcal{A}(U)$ and $t\in \mathcal{B}(U)$, with $U$ an open set in
$X$, the assignment \[(s,t) \longmapsto
\overline{\varphi}_U(s)t\]makes $\mathcal{B}(U)$ into a module over
$\mathcal{A}(U)$. What more is that $\mathcal{B}(U)$ is an algebra
over $\mathcal{A}(U)$; in effect,
\[\overline{\varphi}_U(s+ s^\prime)t= \overline{\varphi}_U(s)t+
\overline{\varphi}_U(s^\prime)t\]and \[\overline{\varphi}_U(s)(t+
t^\prime)= \overline{\varphi}_U(s)t+
\overline{\varphi}_U(s)t^\prime\]for all $s, s^\prime\in
\mathcal{A}(U)$ and $t, t^\prime\in \mathcal{B}(U)$. On the other
hand, since the sheafification functor preserves algebraic
structures, cf. [\cite{mallios}(1.54), p.101], $\mathcal{B}$ can be
viewed as an $\mathcal{A}$-algebra sheaf. The $\mathcal{A}$-algebra
sheaf $\mathcal{B}$ thus obtained is called an
\textbf{$\mathcal{A}$-algebra sheaf with respect to the morphism
$\varphi: \mathcal{A}\longrightarrow \mathcal{B}$}, above. More
accurately, we may yet say that $\mathcal{B}$ is an
\textbf{$\varphi(\mathcal{A})$-algebra sheaf}, with
$\varphi(\mathcal{A})\subseteq \mathcal{B}$, as above. Equivalently,
by a \textbf{$\varphi(\mathcal{A})$-algebra sheaf} (or
\textbf{$\varphi(\mathcal{A})$-algebra}, as a shorthand) , we shall
always mean a \textit{morphism $\mathcal{A}\longrightarrow
\mathcal{B}$ of sheaves of algebras, as above}.

Now, let $\mathcal{A}$ be a sheaf of unital commutative
$\mathbb{C}$-algebras, $\mathcal{E}$ an $\mathcal{A}$-module, and
$\mathcal{R}$ a $\varphi(\mathcal{A})$-algebra. A
\textbf{representation} of $\mathcal{R}$ in $\mathcal{E}$ is an
$\mathcal{A}$-morphism
\[\Theta: \mathcal{R}\rightarrow
\mathcal{E}nd_\mathcal{A}(\mathcal{E}),\]which makes the diagram
\[\xymatrix{\mathcal{R}\ar[rr]^{\hspace{-3mm}\Theta} & &
\mathcal{E}nd_\mathcal{A}(\mathcal{E})\\ &
\mathcal{A}\ar[ul]\ar[ru]}\]commutative; the morphism
$\mathcal{A}\longrightarrow \mathcal{E}nd_\mathcal{A}(\mathcal{E})$
in the above diagram is given by \[a\longmapsto aI|_U=aI_U=
(a_VI_V)_{U\supseteq V, open},\]where $a\in \mathcal{A}(U)$ and
$I_U: \mathcal{E}(U)\longrightarrow \mathcal{E}(U)$ with $U$ open in
$X$, denotes the identity $\mathcal{A}(U)$-morphism. Besides, for
all open set $V$ in $U$, $a_V\in \mathcal{A}(V)$, and $s\in
\mathcal{E}(V)$, \[(a_VI_V)(s)= a_Vs\in \mathcal{E}(V).\]

We observe that for all open $U\subseteq X$, $\mathcal{E}(U)$ may be
viewed as a \textit{module} over
$\mbox{Hom}_{\mathcal{A}|_U}(\mathcal{E}|_U, \mathcal{E}|_U)$.
Indeed, let $f\equiv (f_V)\in
\mbox{Hom}_{\mathcal{A}|_U}(\mathcal{E}|_U, \mathcal{E}|_U)$, where
$V$ runs over the open subsets of $U$, and $s\in \mathcal{E}(U)$.
The action \[(f, s)\longmapsto f_U(s)\] defines a
$\mbox{Hom}_{\mathcal{A}|_U}(\mathcal{E}|_U, \mathcal{E}|_U)$-module
structure on $\mathcal{E}(U)$, as was to be shown. So, by means of
the sheafification process, $\mathcal{E}$ may be viewed as an
$\mathcal{E}nd_\mathcal{A}\mathcal{E}$-\textit{module}. Furthermore,
given a representation $\Theta: \mathcal{R}\longrightarrow
\mathcal{E}nd_\mathcal{A}\mathcal{E}$ of a
$\varphi(\mathcal{A})$-algebra sheaf $\mathcal{R}$ in $\mathcal{E}$,
it turns out that $\mathcal{E}$ may be viewed as an
$\mathcal{R}$-module, with the operation of $\mathcal{R}(U)$ on
$\mathcal{E}(U)$ being given by
\[(s, e)\longmapsto \Theta_U(s)e\equiv \Theta_U(s)(e)\] for $s\in
\mathcal{R}(U)$ and $e\in \mathcal{E}(U)$. Like in \cite{lang},
p.554, we write $se$ instead of the more accurate notation
$\Theta_U(s)(e)$.

\begin{mydf} \emph{Let $\mathcal{B}$ be a \textit{sheaf of $\mathbb{C}$-algebras} over a
topological space $X$, and $\mathcal{A}$ a subsheaf of
$\mathcal{B}$. For all open set $U$ in $X$, let
\[\mathcal{A}(U)[t]\]denote the \textit{ring of polynomials}, in the variable
$t$, whose \textit{coefficients} are the \textit{$($local$)$
sections}  of $\mathcal{A}$ on $U$. Similarly, let
\[\mathcal{A}(U)[s],\]where $s\in \mathcal{B}(U)$, denote the
subring of $\mathcal{B}(U)$ of all \textit{polynomial values}
$p(s)$, with $p\in \mathcal{A}(U)[t]$. A \textit{local section}
$s\in \mathcal{B}(U)$ is called \textbf{transcendental} over
$\mathcal{A}(U)$ if the evaluation map
\[\begin{array}{ll} \mathcal{A}(U)[t]\longrightarrow
\mathcal{A}(U)[s], & p\longmapsto p(s)\end{array}\]is an
\textit{isomorphism}.}\end{mydf}

Now, suppose that $\mathcal{A}$ is a \textit{sheaf of unital
commutative of $\mathbb{C}$-algebras} on the topological space $X$.
As above, we let $\mathcal{A}(U)[t]$, where $U\subseteq X$ is open,
be the polynomial ring. It is clear that the correspondence
\begin{equation}\label{poly} U\longmapsto
\mathcal{A}(U)[t],\end{equation}where $U\subseteq X$ is
\textit{open}, yields a \textit{complete presheaf of
$\mathcal{A}$-modules on $X$}. So, $(\ref{poly})$ defines a
\textit{complete $\Gamma(\mathcal{A})$-presheaf on $X$}. The sheaf
generated, on $X$, by the complete presheaf defined by
$(\ref{poly})$ is called the \textbf{sheaf of germs of polynomials
on $\Gamma(\mathcal{A})$}, and is denoted by $$\mathcal{A}[t],$$
where $t$ is a variable. It is easy to verify that the sheaf
$\mathcal{A}[t]$ of germs of polynomials is an
\textit{$\mathcal{A}$-module} on $X$. Thus, based on
\cite{mallios}, Proposition 11.1, p. 51, one obtains, for every
open set $U\subseteq X$, \[\mathcal{A}[t](U)=
\mathcal{A}(U)[t],\]up to an $\mathcal{A}(U)$-isomorphism.

Let $(X, \mathcal{A})$ be an ordered algebraized space, equipped
with the \textit{inverse-positive-section condition}, $\mathcal{E}$
a vector sheaf of rank $n$, and $\varphi\in
\mbox{E}nd_\mathcal{A}\mathcal{E}$. Let $t$ be transcendental over
$\mathcal{A}(U)$, with $U$ open in $X$, and
\[\mathcal{A}(U)[t]\longrightarrow
\mathcal{A}(U)[\varphi_U]\subseteq
\mbox{Hom}_{\mathcal{A}(U)}(\mathcal{E}(U), \mathcal{E}(U))\]be a
representation of the polynomial ring $\mathcal{A}(U)[t]$ in
$\mathcal{E}(U)$. Like in \cite{lang}, p. 561, we have for every
\textit{open} set $U\subseteq X$, a homomorphism
\[\mathcal{A}(U)[t]\longrightarrow \mathcal{A}(U)[\varphi_U],\]which
is obtained by \textit{substituting} $\varphi_U$ for $t$ in
\textit{polynomials}. The $\mathcal{A}(U)$-algebra
$\mathcal{A}(U)[\varphi_U]$ is the \textit{subalgebra} of
$\mbox{End}_{\mathcal{A}(U)}\mathcal{E}(U)$, generated by
$\varphi_U$, and is \textit{commutative} because
\[\varphi^p_U\circ \varphi^q_U= \varphi_U^q\circ \varphi^p_U,\]for
all $p, q\in \mathbb{N}$. Thus, for all $s\in \mathcal{E}(U)$ and
$f(t)\in \mathcal{A}(U)[t]$, where $U$ is, as usual, an \textit{open
subset} of $X$, we put
\[f(t)s\equiv f(t)(s):= f(\varphi_U)(s)\equiv
f(\varphi_U)s;\]consequently $\mathcal{E}(U)$ turns out to be a
\textit{module} over $\mathcal{A}(U)[t]$. Let $M_U$ be any
\textit{$n\times n$ matrix} in $\mathcal{A}(U)$ (for instance the
matrix representing the $\mathcal{A}(U)$-endomorphism $\varphi_U$
relative to a \textit{canonical basis} $\{e_i^U\}$ of
$\mathcal{E}(U)$, where $U$ is a local gauge of the vector sheaf
$\mathcal{E}$, and $e_i^U= \varepsilon_i^U\circ \varphi^U_i$, with
$\varphi^U$ being the $\mathcal{A}|_U$-isomorphism (\ref{gauge}).
The basis $\{e_i^U\}$ is called a \textit{canonical gauge} of
$\mathcal{E}(U)$.). We define the \textbf{characteristic polynomial
section} $P_{M_U}(t)$ of $M_U$ or of $\varphi_U$ to be the
\textit{determinant}
\[\overline{\partial et}_U(tI_U- M_U):= {\det}_U(tI_U- M_U)\in
\mathcal{A}(U)[t],\]where $I_U$ is the unit $n\times n$-matrix
(here, $1:=1_U$ is the (local) identity section over $U$). (We refer
to \cite{mallios}, pp 294-298, for the sheaf-theoretic notation of
the determinant morphism.) Next, we decree that the
\textbf{characteristic polynomial of an endomorphism} $\varphi\in
\mbox{E}nd_\mathcal{A}\mathcal{E}$ is the endomorphism
$P_\varphi(t)\in \mbox{E}nd_\mathcal{A}\mathcal{A}$, given by
\[P_\varphi(t)\equiv (P_{\varphi_U}(t))= (\det_U(tI_U- M_U)),\]where $M_U$
represents $\varphi_U\in \mbox{E}nd_{\mathcal{A}(U)}\mathcal{E}(U)$
with respect to the local gauge $\{e^U_i\}$ of $\mathcal{E}(U)$. Let
\[M_x= M_U(x);\]its characteristic polynomial is, as obviously
expected, given by
\[\overline{{\partial}et}_U(tI_U- M_U)(x):= \overline{{\partial}et}_x(tI_x- M_x)=
{\det}_x(tI_x- M_x)\in \mathcal{A}_x[t],\]where $\mathcal{A}_x[t]:=
\lim_{\overrightarrow{U\ni x}}\mathcal{A}(U)[t]$.

The characteristic polynomial section will also be referred to
simply as \textit{characteristic polynomial}.

We further look at the following correspondence
\begin{equation}\label{chara}U\longmapsto ChP(\mathcal{A}(U)),\end{equation}where
$U$ ranges over the \textit{open} subsets of $X$, while the range of
$(\ref{chara})$ is the set of all \textit{characteristic
polynomials} $P_{M_U}(t)$ of \textit{$n\times n$-matrices} $M_U$,
whose \textit{entries} are \textit{$($local$)$ sections of
$\mathcal{A}$ on $U$}. Now, the \textit{presheaf of full matrix
$\mathbb{C}$-algebras} on $X$,
\begin{equation}\label{chara1}U\longmapsto
M_n(\mathcal{A}(U)),\end{equation}where $U\subseteq X$ is
\textit{open}, and $M_n(\mathcal{A}(U))$ the \textit{$($full$)$
algebra of $n\times n$-matrices}, with entries the (local) sections
of $\mathcal{A}$ on $U$, is complete. Hence, (\ref{chara}) yields a
complete presheaf on $X$; it is called \textit{presheaf of
characteristic polynomials on $X$}. For the restriction maps of the
presheaf defined by (\ref{chara1}), see \cite{mallios}, pp 280-281.

So, we denote by \[\mathcal{C}h\mathcal{P}(\mathcal{A})\]the
\textit{sheaf of modules on $X$}, generated by the previous
presheaf.

Now before we proceed over to the version of the
\textit{Cayley-Hamilton theorem} in this setting, we signal in
passing (One might work through all the details to their own
satisfaction!) that all the fundamental classical properties of the
\textit{determinant morphism} are also valid in our context. One of
the properties, useful for the scope of the present paper, follows
after this: Let $(X, \mathcal{A}, \mathcal{P})$ be as usual an
ordered $\mathbb{R}$-algebraized space satisfying the
inverse-positive-section condition. Let $A= (s_{ij})\in
M_n(\mathcal{A}(X))$, and $\widetilde{A}= (t_{ij})\in
M_n(\mathcal{A}(X))$ such that \[t_{ij}:=
(-1)^{i+j}\overline{{\partial}et}_X(A_{ji})\equiv
(-1)^{i+j}{\det}_X(A_{ji});\] $A_{ji}$ is the $(n-1)\times (n-1)$
matrix obtained from $A=(s_{ij})_{1\leq i, j\leq n}$ by deleting the
$j$-th row and $i$-th column.

\begin{prop} $($\textbf{Laplace decomposition}$)$, cf. $\cite{remmert}.$ \emph{Let $\det_X(A)= s$, with $s\in \mathcal{A}^\bullet(X)$.
Then, $A\widetilde{A}= \widetilde{A}A= sI$. Furthermore,
$\det_X(A)\in \mathcal{A}^\bullet(X)$ if and only if $A\in
M_n(\mathcal{A})^\bullet(X)$; consequently
\[A^{-1}= s^{-1}\widetilde{A}.\]}\end{prop} \P The proof goes along
similar lines as the proof of Proposition 4.16 in \cite{lang}.
\proofterminator

\begin{theo}\emph{$($\textbf{Cayley-Hamilton}$)$ Let $\mathcal{E}$ be a \textit{vector sheaf of rank
$n$} on $X$, and $\varphi: \mathcal{E}\rightarrow \mathcal{E}$ an
$\mathcal{A}$-morphism. Then,
\[P_\varphi(\varphi)\equiv (P_{\varphi_U}(\varphi_U))=(0_U)\equiv
0.\]}
\end{theo} \P Here as well, we base our proof on the proof of
Theorem 3.1., \cite{lang}, p.561. In fact, let $U$ be an
\textit{open} subset of $X$ such that $\mathcal{E}|_U=
\mathcal{A}^n|_U$, and $\{e_i^U\}_{1\leq i\leq n}$ be the gauge on
$\mathcal{E}|_U(U)= \mathcal{E}(U)$, corresponding  the Kronecker
gauge of . Then, since $\mathcal{E}(U)$ may be viewed as a module
over $\mathcal{A}(U)[t]$, one has \[te_j^U=
\sum_{i=1}^ns_{ij}^Ue_i^U,\]where $1\leq j\leq n$, and
$(s_{ij}^U)\equiv M_U$ is the matrix representing $\varphi_U$ with
respect to $\{e_i^U\}$. Let $B_U(t)= tI_U- M_U$; then
\[\widetilde{B}_U(t)B_U(t)= P_{\varphi_U}(t)I_U,\]and
\[\widetilde{B}_U(t)B_U(t)\left(\begin{array}{c}e_1^U \\
\vdots \\ e_n^U\end{array}\right)= \left(\begin{array}{c}P_{\varphi_U}(t)e_1^U \\
\vdots \\ P_{\varphi_U}(t)e_n^U\end{array}\right)= \left(\begin{array}{c}0_U \\
\vdots \\ 0_U\end{array}\right)\]because \[B(t)\left(\begin{array}{c}e_1^U \\
\vdots \\ e_n^U\end{array}\right)= \left(\begin{array}{c}0_U \\
\vdots \\ 0_U\end{array}\right).\] It follows that
$P_{\varphi_U}(t)(\mathcal{E}(U))= \{0_U\}$, and therefore
$P_{\varphi_U}(\varphi_U)(\mathcal{E}(U))= \{0_U\}$. This implies
that $P_{\varphi_U}(\varphi_U)(\mathcal{E}(U))= \{0_U\}$, as was to
be shown. \proofterminator

Now, suppose that the pair $(X, \mathcal{A})$ is an ordered
algebraized space, satisfying the \textit{inverse-positive-section
condition}. Moreover, as above, we suppose that $\mathcal{E}$ is a
vector sheaf of rank $n$ on $X$, and $\varphi:
\mathcal{E}\longrightarrow \mathcal{E}$ an
$\mathcal{A}$-endomorphism; and we let $s\in \mathcal{E}(U)\equiv
\Gamma(U, \mathcal{E})$, where $U$ is an open set in $X$. We
further consider the \textit{associated} endomorphism
$\overline{\varphi}\equiv (\overline{\varphi}_U)_{U\in
\mathcal{T}}: \Gamma(\mathcal{E})\longrightarrow
\Gamma(\mathcal{E})$, where $\mathcal{T}$ is the assumed topology
on $X$, and $\Gamma(\mathcal{E})$ is the presheaf of sections of
$\mathcal{E}$. Now, we fix $U\in \mathcal{T}$  such that
$\mathcal{E}|_U= \mathcal{A}^n|_U$. By an \textbf{{eigenvector
section}}, or just \textbf{{eigenvector}}, of
$\overline{\varphi}_U\in
\mbox{End}_{\mathcal{A}(U)}(\mathcal{E}(U))$, we mean a
\textit{nowhere-zero (local) section} $s\in \mathcal{E}(U)$, such
that there exists a \textit{section} $\lambda\in \mathcal{A}(U)$
for which we have \begin{equation} \label{evector}
\overline{\varphi}_U(s)= \lambda s,\end{equation} or equivalently
\[M_Us= \lambda s,\]where $M_U$ is the \textit{matrix} representing
$\overline{\varphi}_U$ with respect to the local \textit{Kronecker
gauge} on $\mathcal{E}(U)$. (A (local) section $s\in
\mathcal{E}(U)$ is called a \textit{nowhere-zero} section if
$s_x\equiv s(x)\neq 0_x$ for all $x\in U$.) The \textit{scalar
section} $\lambda$, in Equation (\ref{evector}), is called an
\textbf{eigenvalue section} or simply \textbf{eigenvalue} of the
morphism $\overline{\varphi}_U$.

Next, we construct the following presheaf of sets, namely
\begin{equation}\label{pv}U\longmapsto
\mbox{PV}(\mathcal{E}(U)),\end{equation}where $U$ varies over the
\textit{open subsets} of $X$, and the range of (\ref{pv}), i.e.
$\mbox{PV}(\mathcal{E}(U))$, is the \textit{set of all eigenvectors}
of $\overline{\varphi}_U$, with $\overline{\varphi}_U\in
\mbox{End}_{\mathcal{A}(U)}(\mathcal{E}(U))$.

\begin{prop}\emph{Let $\mathcal{E}$ be an $\mathcal{A}$-module on $X$.
The presheaf $(\mbox{PV}(\mathcal{E}(U)), \sigma^U_V)$, where for
$s\in \mbox{PV}(\mathcal{E}(U))$, $\sigma^U_V(s)\equiv s|_V$, is a
\textit{complete} presheaf.} \end{prop}

\P Indeed, let $\mathcal{U}=\{U_\alpha\}_{\alpha\in I}$ be an open
covering of $U$, and let $s, t$ be two elements of
$\mbox{PV}(\mathcal{E}(U))$ such that
\[\begin{array}{ll}\sigma^U_{U_\alpha}(s)\equiv s_\alpha=t_\alpha\equiv
\sigma^U_{U_\alpha}(t), & \alpha\in I,\end{array}\]where the
$\sigma^U_V$, $U\supseteq V=\mbox{open}$, are the
\textit{restriction maps} of the aforementioned presheaf. Now, as
before let $\Gamma(\mathcal{A})\equiv (\Gamma(U, \mathcal{A}),
\rho^U_V)$ be the presheaf of sections of the sheaf $\mathcal{A}$.
Then, we have, assuming that $\lambda\in \mathcal{A}(U)$ is the
eigenvalue associated with the eigenvector $s\in
\mbox{PV}(\mathcal{E}(U))$, that
\[\overline{\varphi}_{U_\alpha}(s_\alpha)=
\sigma^U_{U_\alpha}(\overline{\varphi}_U(s))=
\sigma^U_{U_\alpha}(\lambda s)=
\rho^U_{U_\alpha}(\lambda)\sigma^U_{U_\alpha}(s)=\rho^U_{U_\alpha}(\lambda)s_\alpha\equiv
\lambda_\alpha s_\alpha.\]Likewise,
\[\overline{\varphi}_{U_\alpha}(t_\alpha)= \mu_\alpha
t_\alpha,\]with $\mu\in \mathcal{A}(U)$ being the eigenvalue for the
eigenvector $t$. As, by \textit{hypothesis}, $s_\alpha= t_\alpha$,
it follows that $\overline{\varphi}_{U_\alpha}(s_\alpha)=
\overline{\varphi}_{U_\alpha}(t_\alpha)$; whence
$\rho^U_{U_\alpha}(\lambda)\equiv\lambda_\alpha= \mu_\alpha\equiv
\rho^U_{U_\alpha}(\mu)$. But $\Gamma(\mathcal{A})$ is a
\textit{complete} presheaf, therefore $\lambda= \mu$; so that $s=t$,
as desired.

Now, let $(s_\alpha)\in
\prod_\alpha\mbox{PV}(\mathcal{E}(U_\alpha))$ such that for any
$U_{\alpha\beta}\equiv U_\alpha\cap U_\beta\neq \emptyset$ in
$\mathcal{U}$, one has
\begin{equation}\sigma^{U_\alpha}_{U_{\alpha\beta}}(s_\alpha)\equiv
s_\alpha|_{U_{\alpha\beta}}= s_\beta|_{U_{\alpha\beta}}\equiv
\sigma^{U_\alpha}_{U_{\alpha\beta}}(s_\beta). \label{ax2}
\end{equation}
The sequence $(s_\alpha)$ of \textit{eigenvectors} gives rise to a
sequence \[\left(M_{n,\alpha}\right)\in \prod_\alpha
M_n(\mathcal{E}(U_\alpha))\] of $n\times n$-matrices whose entries
are (local) \textit{sections} of $\mathcal{E}$, and admitting the
$s_\alpha$ as eigenvectors correspondingly. It is clear that for any
$\alpha$, $\beta\in I$ such that $s_\alpha$, $s_\beta$ fulfill
(\ref{ax2}), one has
\[M_n(\sigma^{U_\alpha}_{U_{\alpha\beta}})(s_{ij}^{U_\alpha}):=
(\sigma^{U_\alpha}_{U_{\alpha\beta}}(s_{ij}^{U_\alpha}))=
(\sigma^{U_\beta}_{U_{\alpha\beta}}(s_{ij}^{U_\beta}))=:
M_n(\sigma^{U_\beta}_{U_{\alpha\beta}})(s_{ij}^{U_\beta}),\]where
$s_\alpha$ and $s_\beta$ are eigenvectors of matrices
$(s_{ij}^{U_\alpha})\in M_n(\mathcal{E}(U_\alpha))$ and
$(s_{ij}^{U_\beta})\in M_n(\mathcal{E}(U_\beta))$, respectively.
But the \textit{presheaf } \[(M_n(\mathcal{E}(U)),
M_n(\sigma^U_V)),\]cf. \cite{mallios}, p. 281, is a
\textit{complete} presheaf, therefore there exists a matrix $M\in
M_n(\mathcal{E}(U))$ such that \[M_n(\sigma^U_{U_\alpha})(M)=
M_{n,\alpha}\]for all $\alpha\in I$. Let $s\in \Gamma(U,
\mathcal{E})$ such that $\sigma^U_{U_\alpha}(s)= s_\alpha$,
$\alpha\in I$. It is easily seen that \[\sigma^U_{U_\alpha}(Ms)=
M_n(\sigma^U_{U_\alpha})(M)\sigma^U_{U_\alpha}(s)=
M_{n,\alpha}s_\alpha= \lambda_\alpha s_\alpha,\] $\alpha\in I$,
which implies that \[Ms= \lambda s,\]where $\lambda\in
\mathcal{A}(U)$ is derived from the $\lambda_\alpha$, $\alpha\in
I$, with $\lambda_\alpha\in \mathcal{A}(U_\alpha)$ and
$\sigma^{U_\alpha}_{U_{\alpha\beta}}(\lambda_\alpha)=
\sigma^{U_\beta}_{U_{\alpha\beta}}(\lambda_\beta)$. Hence, axiom
(S2), see \cite{mallios}, p46, is satisfied.
 \proofterminator

\begin{mydf} \emph{Let $\mathcal{A}$ be a unital commutative
$\mathbb{C}$-algebra sheaf on a topological space $X$, and let
$\mathcal{E}$ be an $\mathcal{A}$-module on $X$. We denote by
\[\mathcal{P}\mathcal{V}(\mathcal{E})\] the sheaf on $X$, generated
by the presheaf defined by $(\ref{pv})$. We call it the
\textbf{eigenvector sheaf} of $\mathcal{E}$ or \textbf{sheaf of
germs of eigenvectors} of $\mathcal{E}$. }\end{mydf}

\begin{prop} \emph{Let $\mathcal{A}$ be a unital $\mathbb{C}$-algebra
sheaf on a topological space $X$, $\omega: \mathcal{A}^{2n}\oplus
\mathcal{A}^{2n}\longrightarrow \mathcal{A}$, $n\in \mathbb{N}$, a
\textit{symplectic structure} on $\mathcal{A}$, and $\varphi\in
\mathcal{S}p\ \mathcal{A}^{2n}(X)$, cf. Lemma \ref{lem21} for the
definition of $\varphi\in \mathcal{S}p\ \mathcal{A}^{2n}(X)$.
Moreover, let $\lambda\in \mathcal{A}^\bullet(X)$ be an
\textit{eigenvalue of $\varphi$}. Then, $\frac{1}{\lambda}\in
\mathcal{A}(X)$ is an \textit{eigenvalue of $\varphi$}
too.}\end{prop}

\P Let $\{s_i\}_{1\leq i\leq 2n}$ be a \textit{basis} of
$\mathcal{A}^{2n}(X)$ such that $(\omega_{ij})_{1\leq i, j\leq
2n}= J$, where $J$ is given by (\ref{matrix}) and $\omega_{ij}=
\omega(s_i, s_j)$, and let $M$ be the $2n\times 2n$-\textit{matrix
representing} the symplectomorphism $\varphi$ with respect to the
aforementioned basis.

Consider the \textit{characteristic polynomial (section)} of $M$
\[P(t)= {\det}_X(M-tI),\]where $I$ is understood as the $2n\times
2n$ identity matrix, and $t\in \mathcal{A}(X)$ a \textit{variable}.
Then, by virtue of (\ref{unit}) and (\ref{cor24}), we have \[P(t)=
t^{2n}P(\frac{1}{t}),\]as is done in \cite{ralph} and \cite{berndt}.
Thus, \[P(\lambda)= \lambda^{2n}P(\frac{1}{\lambda}).\]But
$P(\lambda)=0$ by Cayley-Hamilton theorem, and $\lambda \in
\mathcal{A}(X)= \mathcal{A}^\bullet(X)$, so that
$P(\frac{1}{\lambda})= 0$. \proofterminator

\end{document}